
\input amstex
\documentstyle{amsppt}
\magnification=1200

\hoffset=-0.5pc
\vsize=57.2truepc
\hsize=38truepc
\nologo
\spaceskip=.5em plus.25em minus.20em

\define\cor{\roman{ret}}
\define\rh{r}
\define\wide{}
\define\g{{\frak g}}
\define\Bobb{\Bbb}
 \define\atiyboo{1}

 \define\atibottw{2}

\define\biecktwo{3}
\define\bisguron{4}
\define\bisgurtw{5}
 \define\bottone{6}
\define\botshust{7}
\define\brownboo{8}
\define\eilmac{9}
\define\fockrotw{10}
\define\goldmone{11}
\define\gurupone{12}
\define\gururaja{13}
\define\helgaboo{14}
   \define\modus{15}
\define\modustwo{16}

\define\singula{17}
\define\singulat{18}
\define\smooth{19}
\define\poisson{20}
\define\locpois{21}
    \define\srni{22}
\define\huebjeff{23}
\define\jeffrone{24}
\define\jeffrtwo{25}
\define\jeffrthr{26}
\define\karshone{27}
\define\maclaboo{28}

\define\narasesh{29}

\define\shulmone{30}
\define\sjamlerm{31}
\define\trottone{32}
 \define\weiltwo{33}
\define\weinstwe{34}

\topmatter
\title
Group systems, groupoids, and
moduli spaces of parabolic bundles
\endtitle
\author
K. Guruprasad, J. Huebschmann, L. Jeffrey, and A. Weinstein
\endauthor
\affil
Department of Mathematics
\\
 Indian Institute of Science
\\
Bangalore-560 012, INDIA
\\
kguru\@math.iisc.ernet.in
\\
\phantom{bbb}
\\
Max Planck Institut f\"ur Mathematik
\\
Gottfried Claren Str. 26
\\
D-53 225 BONN
\\
huebschm\@mpim-bonn.mpg.de
\\
\phantom{bbb}
\\
Department of Mathematics
\\
McGill University
\\
Burnside Hall
\\
805 Sherbrooke St. West
\\
Montreal QC Canada H3A 2K6
\\
jeffrey\@gauss.math.mcgill.ca
\\
\phantom{bbb}
\\
Department of Mathematics
\\
University of California
\\
Berkeley, CA 94720-3840 USA
\\
alanw\@math.berkeley.edu
\endaffil
\date{August 3, 1995}
\enddate
\endkeywords
\endsubjclass

\endtopmatter
\document

\beginsection Introduction

Moduli spaces of homomorphisms
or more generally twisted homomorphisms
from fundamental groups of
surfaces to compact connected Lie groups
were connected with geometry
through their
identification with moduli spaces of holomorphic vector bundles
\cite\narasesh.
Atiyah and Bott
\cite\atibottw\
initiated a new approach to the study of these
moduli spaces
by identifying them with moduli spaces
of projectively flat
constant central curvature
connections
on principal bundles over Riemann
surfaces,
which they analyzed by methods of gauge theory.
In particular, they showed that
an invariant inner product on the
Lie algebra
of the Lie group in question
induces
a natural symplectic structure
on a certain smooth open stratum.
Although this moduli space
is a finite dimensional object,
generally a stratified space which is locally semi algebraic
\cite\smooth\
but sometimes a manifold,
its symplectic structure
(on the stratum just mentioned)
was obtained by applying the method of symplectic
reduction to the action of an infinite dimensional group
(the group of gauge transformations)
on an infinite dimensional symplectic manifold
(the space of all connections on a principal bundle).
\smallskip
This infinite-dimensional
approach to moduli spaces has deep roots in quantum field theory
\cite\atiyboo,
but it is nevertheless interesting to try to
avoid the technical difficulties of infinite dimensional analysis
by using purely finite dimensional methods to construct the symplectic
structure and to derive some of its properties.
This also allows for arbitrary, not necessarily compact,
Lie groups.
This program has been carried forward by several authors in the past ten years,
with the result being not only technical simplification,
but also new insight into the geometry of the moduli spaces,
especially into their singularities \cite{\singula\ -- \locpois}.
See \cite\srni\ for a leisurely introduction.
\smallskip
To date, most of the program just described has been worked out
only for compact Riemann surfaces
without boundary;
see however \cite\modustwo.
The purpose of this article is to
extend
these results and methods to the case of Riemann surfaces
with a finite number of punctures or, equivalently,
with a finite number of boundary components,
corresponding to the study of parabolic vector bundles in the
holomorphic category.
Specifically, we deal with the results listed below;
the references
indicate sources for the closed compact case
except \cite\modustwo\ (see below).
\roster
\item"$\bullet$"
A description of the symplectic form in terms of the cup product
on the cohomology of the fundamental group
of the surface
in question
with values in the Lie algebra
\cite\goldmone.
\item"$\bullet$"
A proof, using a double complex of Bott and Shulman rather than gauge theory,
that the form constructed by using group cohomology is closed
\cite\weinstwe,
thereby allowing for a general Lie group, not necessarily compact.
\item"$\bullet$"
A proof, using the Bott-Shulman complex,
that the moduli space can be obtained by symplectic reduction
from a finite-dimensional symplectic manifold
\cite{\modus,\, \modustwo,\, \huebjeff,\, \jeffrtwo}.
\endroster
Some further historical comments may be in order.
Regarding the second item above, a proof that the symplectic
form is closed, using group cohomology rather than gauge theory,
was originally given by Karshon \cite\karshone;
her proof was reformulated in \cite\weinstwe\
in terms of the double complex
of Bott \cite\bottone\  and Shulman \cite\shulmone.
A partly finite dimensional construction
of the moduli space
was accomplished earlier
by Huebschmann \cite\poisson\
and Jeffrey \cite\jeffrone,
but in these papers infinite dimensional techniques
could not completely be avoided.
A purely finite dimensional
construction
(for the closed compact case)
was
announced in \cite\huebjeff\
and given
in \cite{\modus,\, \jeffrtwo}.  In
 \cite\jeffrthr\  the  Bott-Shulman
 construction was used to give  representatives  of all
 generators of the cohomology ring of certain moduli spaces of representations
 that are smooth analogues of the moduli space treated in
 \cite\weinstwe.
\smallskip
Passing to the case of punctured surfaces, we note that the symplectic
structure on
a certain top stratum of
the moduli space
in this case (see Section 9 below for details
about the stratification) was constructed using methods of gauge theory
in \cite\bisguron.
A naive attempt to imitate the methods
used in the closed compact case
\cite{\modus,\huebjeff,\, \jeffrtwo}
seems to fail because the concept of fundamental group
is too weak to handle peripheral structures;
the special
case
where the fundamental group
(of a closed surface)
is replaced by an orbifold fundamental group
--- in the vector bundle case, this
corresponds to parabolic bundles with rational weights ---
has been successfully treated
in \cite\modustwo, though.
\smallskip
The principal
innovation in this paper is to
replace the fundamental group
by two more general concepts
which enable us to overcome
the difficulties with
the peripheral structure in general:
by that of
a {\it group system\/}
\cite\trottone\
and that of a suitable {\it fundamental groupoid\/}.
For our purposes, {\it both\/} notions do
{\it not\/}
serve for equivalent purposes; rather, the two {\it complement\/}
each other.
Group systems
provide the appropriate
concept to handle
the global structure of the moduli space
while the fundamental groupoid turns out to be
a crucial tool for a successful
treatment of the infinitesimal structure.
In
fact,
 a compact orientable topological surface
$\Sigma$
with $n\geq 1$ boundary  circles
$S_1,\dots,S_n$ gives rise
to a group system
$(\pi; \pi_1,\dots, \pi_n)$
(see Section 1 below for details on this notion)
with
$\pi = \pi_1(\Sigma),
\pi_j = \pi_1(S_j) \cong \bold Z$,
and with a chosen generator $z_j$ of each $\pi_j$,
referred to henceforth as a {\it surface group system\/}.
Given
a Lie group $G$, not necessarily compact, and
an $n$-tuple
$\bold C = (C_1,\dots,C_n)$
of conjugacy classes
in $G$,
we denote by $\roman{Hom}(\pi,G)_{\bold C}$
the space of homomorphisms $\chi$ from $\pi$ to $G$
for which
the value $\chi(z_j)$ of each generator
$z_j$
of
lies in $C_j$, for $1 \leq j\leq n$.
Given a nondegenerate invariant
symmetric bilinear form on the Lie algebra $\g$ of $G$,
not necessarily positive definite,
we shall construct an {\it extended moduli space\/}
$\Cal M(\wide {\Cal P},G)_{\bold C}$,
that is to say, a smooth
symplectic manifold and
a hamiltonian
$G$-action and momentum mapping
whose reduced space
$\Cal M(\wide {\Cal P},G)_{\bold C}\big / \big / G$
is homeomorphic
to the space
$\roman{Rep}(\pi,G)_{\bold C}$
of representations, the orbit space
for the action of $G$ by conjugation
on $\roman{Hom}(\pi,G)_{\bold C}$.
The global
construction of the space
$\Cal M(\wide {\Cal P},G)_{\bold C}$,
of its closed
2-form,
and of the momentum mapping,
involve
the surface group system,
whereas nondegeneracy of the form
is
proved
by relating
the infinitesimal
part of
structure
with Poincar\'e duality
in relative
or more precisely
{\it parabolic\/}
cohomology
of
a certain
fundamental groupoid.
\smallskip
Another proof of nondegeneracy
is given in a companion paper to this one
\cite\gururaja, which uses a metric on the space
of parabolic cocycles, leading to a (new?) metric
on the moduli space.
\smallskip
For compact $G$, in the gauge theory setting,
the tangent space of an arbitrary point of the top stratum
of the moduli space mentioned above
in the case of a punctured surface
can be identified with the image
of compactly
supported de Rham cohomology
in the usual de Rham cohomology with coefficients
in the adjoint bundle,
calculated with reference to the operator $d_A$ of covariant derivative
with respect to a flat connection $A$ representing
the point in question
\cite\bisguron.
The cohomology of \lq\lq group systems\rq\rq\
is the analogue of compactly supported cohomology in the algebraic
setting, and the tangent space
at a point of the top stratum
can be identified with the image of the corresponding
group systems cohomology in the usual cohomology.
In the last section of this paper, we make this explicit
and show that the symplectic structure
on the
top stratum
obtained here by algebraic methods is equivalent
to the construction via gauge theory.
\smallskip
Another
finite-dimensionalization
of the space of flat connections,
inspired by lattice gauge theory,
was introduced by Fock and Rosly
\cite\fockrotw.
\smallskip
The four authors of this paper would like to thank the
Soci\'et\'e Math\'ematique de France
and the CNRS URA 751 at Lille University
for making it possible, through a meeting at Luminy
in March 1994,
for them to be in the same place at the same time
in order to crystallize
their ideas for a joint paper
which was previously conceived
through separate
meetings and correspondence.
In addition, Guruprasad would like to thank the Ecole Polytechnique,
Palaiseau, and the Centre Emile Borel
of the Institut Henri Poincar\'e
for their hospitality
during his visit when a substantial part of this work
was done.
He is in particular grateful to
Prof. Yvette Kosmann-Schwarzbach and Prof. J.P.Bourguignon
for their hospitality and support.
He is also especially thankful to the University of
California, Berkeley,
for providing an excellent atmosphere
for research.
Huebschmann carried out a substantial part
of this work in the framework of the
VBAC research group of EUROPROJ,
first as a member of the CNRS URA 751 at Lille University
and thereafter
during a stay at the Max Planck Institut
f\"ur Mathematik at Bonn.
He wishes to express his gratitude to it and to its director,
Prof. F.\,Hirzebruch,
for hospitality and support.
Jeffrey's work is partially supported
by NSF grant DMS-9306029.
Weinstein would like to thank the Indian Institute
of Science (Bangalore)
for an invitation
which initiated his collaboration
with Guruprasad.
The work begun there
was continued in the excellent circumstances
provided by the Centre
Emile Borel of the
Institut Henri Poincar\'e.
Finally, his work, and a visit
to Berkeley by Guruprasad
for the completion of the manuscript,
was supported by the National Science Foundation Grant
DMS-93-09653.

\beginsection 1. Group systems

Recall that a {\it group system\/}
$(\pi; \phi_1,\dots,\phi_n, \pi_1,\dots, \pi_n)$
consists of a group $\pi$ together with a family
of groups $\pi_j$ and homomorphisms
$\phi_j$ from $\pi_j$ to $\pi$
\cite\trottone.
We shall occasionally refer to the $\pi_j$ as
{\it peripheral\/} groups.
For any such system,
there is a pair of spaces
$(X,\cup Y_j)$
such that $X$ and the $Y_j$ are aspherical,
$\pi = \pi_1(X),
\pi_j =\pi_1(X_j)$, and the maps $\phi_j$
are induced by inclusion.
The {\it (co)homology of the group system\/}
is that of the pair
$(X,\cup Y_j)$.
{\smc Trotter\/} has given a purely algebraic construction \cite\trottone.
To introduce notation we reproduce it briefly:
\smallskip
Let $R$ be an arbitrary commutative ring, taken henceforth as ground ring.
A {\it resolution over a system\/}
$(\pi; \phi_1,\dots,\phi_n, \pi_1,\dots, \pi_n)$
is a pair of $R\pi$-complexes
$(A,B)$ such that
\roster
\item
$A$ is a resolution over $\pi$;
\item
$B$ is the direct sum of complexes $B_j = R\pi \otimes _{R\pi_j}A_j$
where $A_j$ is a resolution over $\pi_j$ and
$R\pi$ is considered a right $R\pi_j$-module via the map $\phi_j$;
\item
$B$ is a $R\pi$-direct summand of $A$.
\endroster
The $A_j$ are referred to as {\it auxiliary resolutions\/};
occasionally we shall refer to $B$ as the {\it peripheral part\/}
of the resolution.
Given
a resolution
$(A,B)$
over a system
$(\pi; \phi_1,\dots,\phi_n, \pi_1,\dots, \pi_n)$, the exact sequence
$$
0
@>>>
B
@>>>
A
@>>>
A/B
@>>>
0
\tag 1.1
$$
of $R\pi$-modules
splits,
and the homology
$\roman H_*(\{\phi_j\},\cdot)$
and cohomology
$\roman H^*(\{\phi_j\},\cdot)$
of the system are defined from the
$R\pi$-complex
$A/B$, which plays the role of a \lq\lq relative\rq\rq\ resolution.
In particular, for every
$R\pi$-module $M$,
the exact sequence (1.1)
gives rise to natural long exact  sequences in homology and cohomology
of the kind
$$
\dots
@>>>
\roman H_k(\{\pi_j\},M)
@>>>
\roman H_k(\pi,M)
@>>>
\roman H_k(\{\phi_j\},M)
@>>>
\roman H_{k-1}(\{\pi_j\},M)
@>>>
\dots
\tag1.2
$$
and
$$
\dots
@<<<
\roman H^{k+1}(\{\phi_j\},M)
@<<<
\roman H^k(\{\pi_j\},M)
@<<<
\roman H^k(\pi,M)
@<<<
\roman H^k(\{\phi_j\},M)
@<<<
\dots  \ .
\tag1.3
$$
Here
$\roman H_k(\{\phi_j\},M)$
and
$\roman H^k(\{\phi_j\},M)$
are just the
direct sums of the homology
and
cohomology groups
$\roman H_k(\pi_j,M)$
and $\roman H^k(\pi_j,M)$,
respectively.
\smallskip
Given a group system
$(\pi; \phi_1,\dots,\phi_n, \pi_1,\dots, \pi_n)$
where the $\phi_j$ are inclusions of subgroups
we shall henceforth
suppress $\phi_1,\dots,\phi_n$
in notation and simply
 write
$(\pi; \pi_1,\dots, \pi_n)$ and, likewise,
we shall write
$
\roman H_*(\pi,\{\pi_j\};\cdot)
$
and
$\roman H^*(\pi,\{\pi_j\};\cdot)$
for
$\roman H_*(\{\phi_j\},\cdot)$
and
$\roman H^*(\{\phi_j\},\cdot)$,
respectively.
\smallskip
The explicit construction of a resolution of a system
can concisely be handled by means of a corresponding
{\it fundamental groupoid\/}.
We explain this for surface group systems in the next section.

\beginsection 2. Surface group systems

Let $\Sigma$ be a compact orientable topological surface
with boundary $\partial \Sigma$ consisting of $n$ circles
$S_1,\dots,S_n$;
we suppose that
when the genus $\ell$ of $\Sigma$
is zero
there are
$n \geq 3$
boundary circles.
This surface gives rise
to a group system
$(\pi; \pi_1,\dots, \pi_n)$
with
$\pi = \pi_1(\Sigma),
\pi_j = \pi_1(S_j) \cong \bold Z$,
referred to henceforth as a {\it surface group system\/}.
When the boundary $\partial \Sigma$
is non-empty
the group $\pi$ is free non-abelian; yet it is convenient to use the
presentation
$$
\Cal P = \langle x_1,y_1,  \dots,x_\ell,y_\ell,z_1,\dots,z_n; r
\rangle,
\tag2.1
$$
where
$$
r = \Pi [x_j,y_j] z_1 \dots z_n.
$$
The Reidemeister-Fox calculus, applied to the presentation $\Cal P$,
yields
the free resolution
$$
\bold R(\Cal P)\colon
R_2(\Cal P)
@>\partial_2>>
R_1(\Cal P)
@>\partial_1>>
R\pi
\tag2.2
$$
of $R$ in the category of left
$R\pi$-modules.
Here
$$
R_2(\Cal P) =R\pi[r], \quad
R_1(\Cal P) =
R\pi[x_1,y_1,\dots,x_\ell,y_\ell,z_1,\dots,z_n],
$$
and
the boundary operators
$\partial_j$ are given by the
formulas
$$
\aligned
\partial_1[x_i] &= (x_i-1)
\\
\partial_1[y_i] &= (y_i-1)
\\
\partial_1[z_j] &= (z_j-1)
\\
\partial_2[r] &=
\sum
\frac{\partial r}{\partial x_i} [x_i]
+
\sum
\frac{\partial r}{\partial y_i} [y_i]
+
\sum
\frac{\partial r}{\partial z_j} [z_j]  .
\endaligned
$$
Here and henceforth
we decorate the free generators
of the modules coming into play in the resolution
by square brackets,
to distinguish them
from the corresponding elements
of the group $\pi$ etc.
The chain complex
arising from (2.2)
which calculates the absolute homology of
$\pi$ with values in $R$ comes down to
$$
R[r]
@>\overline \partial_2>>
R[x_1,y_1,\dots,x_\ell,y_\ell,z_1,\dots,z_n]
@>{0}>>
R,\quad
\overline \partial_2[r] = [z_1] + \dots + [z_n].
\tag2.3
$$
\smallskip
A resolution of the group system
is concisely handled by means of
the following groupoid
the full force of which
will be exploited only
in Section 8, though.
Pick a base point $p_0$
not on the boundary
and, moreover,
for each boundary component $S_j$,
pick a base point
$p_j$.
This determines
the  subgroupoid
$\widetilde \pi= \Pi(\Sigma;p_0,p_1,\dots,p_n)$
of the fundamental
groupoid of $\Sigma$ consisting of homotopy classes
of paths in $\Sigma$ with endpoints contained in the set
$\{p_0,p_1,\dots,p_n\}$.
To obtain a presentation of it we decompose
$\Sigma$ into cells as follows,
where we do not distinguish in notation between the chosen
edge paths and their homotopy classes relative
to their end points:
Let $x_1,y_1,\dots,x_{\ell}, y_{\ell}$
be closed paths
which (i) do not meet the boundary,
(ii) have $p_0$ as starting point, and
(iii) yield the generators
respectively
$x_1,y_1,  \dots,x_\ell,y_\ell$
of the fundamental group
$\pi=\pi_1(\Sigma,p_0)$;
for $j = 1,\dots, n$,
let $a_j$
be
the boundary path of the $j$'th
boundary circle,
having $p_j$ as starting point,
and let
$\gamma_j$
be a path
from $p_0$ to $p_j$.
When we cut $\Sigma$ along these
1-cells we obtain a disk $D$ whose boundary
yields the defining relation of
$\widetilde \pi= \Pi(\Sigma;p_0,p_1,\dots,p_n)$.
The resulting
presentation
of
$\widetilde \pi$
looks like
$$
\widetilde {\Cal P} = \langle x_1,y_1,  \dots,x_\ell,y_\ell,a_1,\dots,a_n,
\gamma_1,\dots,\gamma_n; \widetilde r\rangle,
\tag2.4
$$
where
the relator  now reads
$$
\widetilde r = \Pi [x_j,y_j] \Pi \gamma_j a_j \gamma_j^{-1}.
$$
This is consistent with the presentation (2.1) of
the fundamental group if we identify each generator
$z_j$ with $\gamma_j a_j \gamma_j^{-1}$.
Reading the boundary $\widetilde r$ counterclockwise around $D$
determines an orientation of $D$ and hence of
$\Sigma$ in the usual way.
\smallskip
The Reidemeister-Fox calculus, applied to
$\widetilde {\Cal P}$,
yields
the free resolution of $R$
$$
\bold R(\widetilde {\Cal P})\colon
R_2(\widetilde {\Cal P})
@>\partial_2>>
R_1(\widetilde {\Cal P})
@>\partial_1>>
R_0(\widetilde {\Cal P}) = R\pi[p_0,\dots,p_n]
\tag2.5
$$
in the category of left
$R\pi$-modules.
Here
$$
R_2(\widetilde {\Cal P}) =R\pi[\widetilde r], \quad
R_1(\widetilde {\Cal P}) =
R\pi[x_1,y_1,\dots,x_\ell,y_\ell,a_1,\dots,a_n,\gamma_1,\dots,\gamma_n],
$$
and
$$
\aligned
\partial_1[x_i] &= (x_i-1)[p_0],
\\
\partial_1[y_i] &= (y_i-1)[p_0],
\\
\partial_1[a_j] &= (z_j-1)[p_j],
\\
\partial_1[\gamma_j] &= [p_j] - [p_0],
\\
\partial_2[\widetilde r] &=
\sum
\frac{\partial r}{\partial x_i} [x_i]
+
\sum
\frac{\partial r}{\partial y_i} [y_i]
+
\sum
\frac{\partial r}{\partial z_j} [a_j]
+
\sum
\frac{\partial r}{\partial z_j} (1-z_j)[\gamma_j] .
\endaligned
\tag2.6
$$
Notice that this amounts to a concise description
of the cellular chains of the universal cover
$\widetilde \Sigma$ of $\Sigma$
whence it is manifestly
a free resolution.
This description of the chains of the universal cover
will be exploited in Section 10 below.
Alternatively,
observe
that dividing out the contractible $R\pi$-subcomplex
generated
by the $[\gamma_j]$ and $[p_j]-[p_0]$
transforms (2.5) into the free resolution
(2.2) above
The resulting chain complex calculating the homology of
$\pi$
with values in $R$ amounts to
$$
R[r]
@>\overline \partial_2>>
R[x_1,y_1,\dots,x_\ell,y_\ell,a_1,\dots,a_n,\gamma_1,\dots,\gamma_n]
@>{\overline \partial_1}>>
R
\tag2.7
$$
where
$$
\alignat 2
\overline \partial_2 [\widetilde r] &=[a_1] + \dots + [a_n],  &&
\\
\overline \partial_1 [x_i] &= \overline \partial_1 [y_i]= 0,
&& \quad 1 \leq i \leq \ell,
\\
\overline \partial_1 [a_j] &= 0,
&& \quad 1 \leq j \leq n.
\\
\overline \partial_1 [\gamma_j] &= [p_j] - [p_0],
&& \quad 1 \leq j \leq n.
\endalignat
$$
It equals that of cellular chains of $\Sigma$,
when we identify the disk $D$ with $[\widetilde r]$.
Notice that when $n \geq 1$, the orientation of $D$
is determined by the boundary relation
$$
\overline \partial D =[a_1] + \dots + [a_n].
$$
\smallskip
To see that (2.5) yields
a resolution
over our group system,
for
$j = 1,\dots,n$,
let
$A_j$
be the small resolution of $R$ over
the free cyclic group $\pi_j$
determined by the choice of generator $a_j$,
and let
$B_j = R\pi \otimes _{R\pi_j}A_j$; explicitly:
$$
B_j
\colon
R\pi
[a_j]
@>{\partial}>>
R\pi[p_j],\quad
\partial [a_j] = (z_j -1) [p_j].
$$
The resolution (2.5)
plainly contains
the
$R\pi$-complex $B=\oplus_{j=1,\dots,n} B_j$
as a direct summand.
Hence
$(\bold R(\widetilde {\Cal P}),B)$ is a resolution
over our group system.
Notice that for $\ell=0$ and $n=1$ the construction does {\it not\/}
yield a resolution over the corresponding group system.
\smallskip
By construction, the quotient complex
$\bold R(\widetilde {\Cal P},\{\pi_j\})$
calculating the (co)homology of our surface system
arises from
(2.5) by dividing out the subcomplex
generated by the $[a_j]$ and $[p_j]$, for $1 \leq j \leq n$.
Thus it
looks like
$$
{\bold R}(\widetilde {\Cal P},\{\pi_j\})\colon
R\pi[\widetilde r]
@>\partial_2>>
R\pi[x_1,y_1,\dots,x_\ell,y_\ell,\gamma_1,\dots,\gamma_n]
@>\partial_1>>
R\pi ;
\tag2.8
$$
its
boundary operators
$\partial_j$ are given by the
formulas
$$
\aligned
\partial_1[x_i] &= (x_i-1)
\\
\partial_1[y_i] &= (y_i-1)
\\
\partial_1[\gamma_j] &= 1
\\
\partial_2[\widetilde r] &=
\sum
\frac{\partial r}{\partial x_i} [x_i]
+
\sum
\frac{\partial r}{\partial y_i} [y_i]
+
\sum
\frac{\partial r}{\partial z_j} (1-z_j)[\gamma_j] .
\endaligned
\tag2.9
$$
In particular, the chain complex calculating the homology of
the system
with values in $R$ amounts to
$$
R[r]
@>\overline \partial_2>>
R[x_1,y_1,\dots,x_\ell,y_\ell,\gamma_1,\dots,\gamma_n]
@>{\overline \partial_1}>>
R
\tag2.10
$$
where
$$
\alignat 2
\overline \partial_2 [\widetilde r] &=0  &&
\\
\overline \partial_1 [x_i] &= \overline \partial_1 [y_i]= 0,
&& \quad 1 \leq i \leq \ell,
\\
\overline \partial_1 [\gamma_j] &= 1,
&& \quad 1 \leq j \leq n.
\endalignat
$$
Thus the 2-chain
$$
b =\widetilde r
\tag2.11
$$
is a relative 2-cycle,
and
$\roman H_2(\pi,\{\pi_j\};R)$
is isomorphic to $R$, generated
by the class $\kappa$ of $b$.
\smallskip
A comparison map
from
$\bold R(\Cal P)$ to
$\bold R(\widetilde {\Cal P},\{\pi_j\})$
inducing the cohomology map
from
$
\roman H^*(\pi,\{\pi_j\};\cdot)
$
to
$
\roman H^*(\pi,\cdot)
$
occurring in (1.3) above
is given by
the chain map
$$
\CD
\bold R(\Cal P)\colon
@.
R\pi[r]
@>\partial_2>>
R\pi[x_1,y_1,\dots,x_\ell,y_\ell,z_1,\dots,z_n]
@>\partial_1>>
R\pi
\\
@.
@VVV
@VVV
@VVV
\\
{\bold R}(\widetilde {\Cal P},\{\pi_j\})\colon \quad
@.
R\pi[\widetilde r]
@>\partial_2>>
R\pi[x_1,y_1,\dots,x_\ell,y_\ell,\gamma_1,\dots,\gamma_n]
@>\partial_1>>
R\pi
\endCD
\tag2.12
$$
which identifies the elements denoted by the same symbols
in the top and bottom row and sends
$[r]$ to $[\widetilde r]$ and
$[z_j]$ to $(z_j-1)[\gamma_j]$, for $ 1 \leq j \leq n$.
In particular, under the induced chain map from
(2.3) to (2.10),
the boundary value
$[z_1] + \dots + [z_n]$
of
$[r]$ in (2.3)
goes to zero.

\beginsection 3. Poincar\'e duality for surface group systems

The
surface group system
$(\pi; \pi_1,\dots, \pi_n)$
is a
{\it two-dimensional Poincar\'e duality
system\/}
over $R$,
that is, a $\roman{PD}^2$-{\it pair\/} in the terminology of
{\smc Bieri-Eckmann} \cite\biecktwo,
having fundamental class
$\kappa \in \roman H_2(\pi,\{\pi_j\};R)$,
so that,
for every $R \pi$-module $M$,
cap product with $\kappa$
yields natural isomorphisms
$$
\cap \kappa
\colon
\roman H^*(\pi,M)
@>>>
\roman H_{2-*}(\pi,\{\pi_j\};M),
\quad
\cap \kappa
\colon
\roman H^*(\pi,\{\pi_j\};M)
@>>>
\roman H_{2-*}(\pi,M),
\tag3.1
$$
cf. \cite\biecktwo,
and these in fact fit into a commutative diagram
$$
\CD
@>>>
\roman H^*(\pi,\{\pi_j\};M)
@>>>
\roman H^*(\pi,M)
@>>>
\roman H^*(\{\pi_j\},M)
@>>>
\\
@.
@V{\cap \kappa}VV
@V{\cap \kappa}VV
@V{\cap \partial\kappa}VV
@.
\\
@>>>
\roman H_{2-*}(\pi,M)
@>>>
\roman H_{2-*}(\pi,\{\pi_j\};M)
@>>>
\roman H_{1-*}(\{\pi_j\},M)
@>>>
\endCD
$$
whose horizontal sequences are the corresponding
long exact homology and cohomology sequences
of the group system.
In particular,
$\roman H^2(\pi,\{\pi_j\};R)$
is also just a copy of $R$.
Geometrically this duality is exactly that of the
surface $\Sigma$ with boundary $\partial \Sigma$ with coefficients
determined by $M$, viewed as a local system, that is,
$$
\cap e
\colon
\roman H^*(\Sigma,M)
@>>>
\roman H_{2-*}(\Sigma,\partial \Sigma;M),
\quad
\cap e
\colon
\roman H^*(\Sigma,\partial \Sigma;M)
@>>>
\roman H_{2-*}(\Sigma,M),
\tag3.2
$$
where
$e \in \roman H_{2}(\Sigma,\partial \Sigma;R)$
refers to the orientation class.
\smallskip
Let
$R=\Bobb R$, the reals, let
$V$ be a finite dimensional real vector space
with a symmetric bilinear form
$\cdot$\,,
and suppose $V$ endowed with a structure
of $\Bobb R\pi$-module preserving the given
symmetric bilinear form
(i.e. the action of $\pi$ preserves the form).
Via the multiplicative structure
of the cohomology of a group system
--- this amounts of course
to the multiplicative structure
of the cohomology of the pair
$(\Sigma,\partial \Sigma)$ ---
the symmetric bilinear form $\cdot$ induces a pairing
$$
\roman H^1(\pi,\{\pi_j\};V)
\otimes
\roman H^1(\pi,V)
@>>>
\roman H^2(\pi,\{\pi_j\};\Bobb R)
\tag3.3
$$
which, combined with
$$
\cap \kappa
\colon
\roman H^2(\pi,\{\pi_j\};\Bobb R)
@>>>
\roman H_0(\pi,\Bobb R) = \Bobb R,
$$
yields a bilinear pairing
$$
\roman H^1(\pi,\{\pi_j\};V)
\otimes
\roman H^1(\pi,V)
@>>>
\Bobb R;
\tag3.4
$$
Poincar\'e duality
in the cohomology of the group system
implies that
the pairing (3.4) is nondegenerate
provided that $\cdot$ is nondegenerate.
The multiplicative structure
has been made explicit in \cite\trottone \
for systems with a single peripheral
subgroup. By means of an appropriate groupoid
we shall make explicit the multiplicative
structure in Section 8 below
in the general case.
\smallskip
Next we write
$\roman H_{\roman{par}}^1(\pi,\{\pi_j\};V)$
for the image of
$\roman H^1(\pi,\{\pi_j\};V)$
in
$\roman H^1(\pi,V)$
under the canonical map, cf. (1.3),
and we refer to it as
(first) {\it parabolic cohomology\/},
with values in $V$.
Parabolic cohomology classes are represented
by
{\it parabolic\/} 1-cocycles, that is,
by
1-cocycles
$\zeta \colon \pi \to V$
having the property that,
for every $z_j,\, 1 \leq j \leq n,$ there is
an element $X_j$ in $V$ such that
$$
\zeta(z_j) = z_j X_j - X_j.
\tag3.5
$$
We denote
the space of
parabolic 1-cocycles
by $Z_{\roman{par}}^1(\pi,\{\pi_j\};V)$.
Parabolic
1-cocycles
and parabolic cohomology
have been introduced by {\smc A. Weil} \cite\weiltwo,
for arbitrary finitely generated planar discontinuous groups,
and he noticed that,
for such a group
with only elliptic
and hyperbolic generators
(which, in the present description,
amounts to
imposing the additional relations
saying that every generator of the kind
$z_j$
has finite order)
every
1-cocycle is parabolic
since the cohomology of a finite group
with coefficients in a real vector space
is trivial.
\smallskip
Using the top row of the commutative
diagram after
(3.1), we have
an exact sequence
$$
0
@>>>
\roman{Ker}(j)
@>>>
\roman H^1(\pi,\{\pi_j\};V)
@>j>>
\roman H^1(\pi,V)
@>>>
\roman{Coker}(j)
@>>>
0
$$
and
the restriction of (3.4) to
$\roman{Ker}(j) \otimes \roman{Im}(j)$ is zero
where
$\roman{Im}(j)$
refers to the image
of $j$ in
$\roman H^1(\pi,V)$.
This implies that the pairing (3.4) yields a pairing
$$
\left(\roman H^1(\pi,\{\pi_j\};V)
\big /\roman{Ker}(j)\right)
\otimes
\roman{Im}(j)
@>>>
\Bobb R.
$$
However
$j$ induces an isomorphism
from
$\roman H^1(\pi,\{\pi_j\};V)\big/\roman{Ker}(j)$
onto
$\roman{Im}(j)$ which equals $\roman H_{\roman{par}}^1(\pi,\{\pi_j\};V)$.
Hence the pairing (3.4)
induces a
skew-symmetric
bilinear pairing
$$
\omega_V
\colon\roman H_{\roman{par}}^1(\pi,\{\pi_j\};V)
\otimes
\roman H_{\roman{par}}^1(\pi,\{\pi_j\};V)
@>>>
\Bobb R.
\tag3.6
$$
When $\cdot$ is nondegenerate,
so is (3.4),
and
$\roman{Im}(j)$
equals the annihilator
in
$\roman H^1(\pi,V)$
of
$\roman{Ker}(j)$.
Hence
(3.6)
is nondegenerate,
that is, a symplectic structure
on
the vector space
$\roman H_{\roman{par}}^1(\pi,\{\pi_j\};V)$,
provided that $\cdot$ is nondegenerate.
This   pairing
will be explicitly calculated
in  Lemma 8.4 below.

\beginsection 4. Representation spaces of surface group systems

Write $F$ for the free group
on the generators in $\Cal P$.
Let $G$ be a Lie group, not necessarily compact,
write $\g$ for its Lie algebra, and let
$\bold C = \{C_1,\dots,C_n\}$ be an $n$-tuple of conjugacy classes
in $G$.
\smallskip
Let $\phi \in \roman{Hom}(F,G)$,
and suppose that
$\phi(r)$ lies in the centre of $G$.
Then the composite of $\phi$
with the adjoint representation of $G$
induces a structure of a (left) $\pi$-module on $\g$,
and we write
$\g_{\phi}$ for $\g$, viewed as as $\pi$-module in
this way.
We shall continue
to take as ground ring $R$ the reals $\Bobb R$.
Application of the functor
$\roman{Hom}_{\Bobb R\pi}(\cdot, \g_{\phi})$
to the free resolution
$\bold R(\Cal P)$
yields the chain complex
$$
\bold C(\Cal P, \g_{\phi})\colon
\roman C^0(\Cal P, \g_{\phi})
@>{\delta_{\phi}^0}>>
\roman C^1(\Cal P, \g_{\phi})
@>{\delta_{\phi}^1}>>
\roman C^2(\Cal P, \g_{\phi}),
\tag4.1
$$
cf. \cite\modus\ (4.1),
computing the group cohomology $\roman H^*(\pi,\g_{\phi})$;
we recall that
there are canonical isomorphisms
$$
\roman C^0(\Cal P, \g_{\phi}) \cong \g,
\quad
\roman C^1(\Cal P, \g_{\phi}) \cong \g^{2\ell+ n},
\quad
\roman C^2(\Cal P, \g_{\phi}) \cong \g .
$$
To recall the geometric significance
of this chain complex,
denote by $\alpha_\phi$
the smooth map
from $G$  to $\roman{Hom}(F,G)$
which assigns $x \phi x^{-1}$ to $x \in G$,
write $\rh$
for the smooth map
from
$\roman{Hom}(F,G)$
to
$G$
induced by the relator $r$
so that the pre-image of the neutral element
$e$ of $G$ equals
the space
$\roman{Hom}(\pi,G)$,
and write
$R_\phi\colon \g^{2\ell+n} \to  \roman T_{\phi} \roman{Hom}(F,G)$
and
$R_{\rh(\phi)}\colon \g \to  \roman T_{\rh (\phi)}G$
for the corresponding operations of right translation.
The tangent maps
$\roman T_e\alpha_{\phi}$
and $\roman T_{\phi}\rh$
make commutative the diagram
$$
\CD
\roman T_eG
@>\roman T_e\alpha_{\phi}>>
\roman T_{\phi} \roman{Hom}(F,G)
@>{\roman T_{\phi} \rh}>>
\roman T_{\rh (\phi)}G
\\
@A{\roman{Id}}AA
@A{\roman R_{\phi}}AA
@A{\roman R_{\rh(\phi)}}AA
\\
\g
@>>{\delta^0_{\phi}}>
\g^{2\ell +n}
@>>{\delta^1_{\phi}}>
\g,
\endCD
\tag4.2
$$
cf. \cite\modus\ (4.2).
The commutativity of this diagram
shows at once that
right translation
identifies the
kernel of
the derivative
$\roman T_{\phi} \rh$
with the kernel
of the coboundary operator
$\delta_{\phi}^1$
from
$\roman C^1(\Cal P, \g_{\phi})$
to
$\roman C^2(\Cal P, \g_{\phi})$,
that is, with the
vector space
$\roman Z^1(\pi,\g_{\phi})$
of $\g_{\phi}$-valued 1-cocycles of $\pi$;
this space does {\it not\/}
depend on a specific presentation $\Cal P$, whence the notation.
We note that $\roman C^1(\Cal P, \g_{\phi}) =
Z^1(F,\g_{\phi})$,
the space of $\g_{\phi}$-valued 1-cocycles for $F$.
\smallskip
For each $j$,
$1 \leq j \leq n$,
write $F_j$
for the subgroup of $F$ generated by $z_j$.
We then have two group systems
$(F; F_1,\dots, F_n)$
and
$(\pi; \pi_1,\dots, \pi_n)$,
together with the obvious morphism
of group systems
from
the former to the latter.
Notice that, for each $j$,
the corresponding homomorphism from
$F_j$ to $\pi_j$
is an isomorphism
but we prefer to maintain a distinction
in notation
between $F_j$ and $\pi_j$.
Extending notation introduced earlier,
we denote by $\roman{Hom}(F,G)_{\bold C}$
the space of homomorphisms $\phi$  from $F$ to $G$
for which
the value $\phi(z_j)$ of each generator $z_j$
lies in $C_j$, for $1 \leq j\leq n$.
The choice of generators induces a decomposition
$$
\roman{Hom}(F,G)_{\bold C}
\cong
G^{2\ell} \times
C_1\times
\dots
\times
C_n.
$$
Abusing notation, we denote
the restriction of $\rh$ to
$\roman{Hom}(F,G)_{\bold C}$
by $\rh$ as well,
so that the pre-image
$\rh^{-1}(e) \subseteq \roman{Hom}(F,G)_{\bold C}$
of the neutral element
$e$ of $G$ equals
the space
$\roman{Hom}(\pi,G)_{\bold C}$.
\smallskip
We now suppose that   our chosen
$\phi \in \roman{Hom}(F,G)$ lies in
$\roman{Hom}(F,G)_{\bold C}$,
viewed as a subspace
of
$\roman{Hom}(F,G)$.
For $j = 1, \dots, n$,
denote by $h_j$ the image in $\g$ of the linear endomorphism
given by
$\roman{Ad}(\phi(z_j)) -\roman{Id}$, so that there results the exact sequence
$$
0
@>>>
\g_j
@>>>
\g
@>>>
h_j
@>>>
0
\tag4.3
$$
of vector spaces,
where
$\g_j$ denotes the
Lie algebra of the
stabilizer
of $\phi(z_j)$;
notice that $h_j$ is the tangent space
of the conjugacy class $C_j$.
\smallskip
The following two observations
will be crucial.

\proclaim{Proposition 4.4}
The values of the operator
$\delta^0_{\phi}$
in {\rm (4.2)}
lie in
$\g^{2\ell} \times h_1 \times \dots\times h_n$,
viewed as a subspace
of
$\roman C^1(\Cal P, \g_{\phi}) \cong \g^{2\ell} \times \g^n$,
and
the first cohomology group of the resulting complex
$$
\bold C_{\roman{par}}(\Cal P, \g_{\phi})\colon
\g
@>{\delta^0_{\phi}}>>
\g^{2\ell} \times h_1 \times \dots\times h_n
@>{\delta^1_{\phi}}>>
\g
\tag4.4.1
$$
equals $\roman H_{\roman{par}}^1(\pi,\{\pi_j\};\g_{\phi})$ .
\endproclaim

\demo{Proof}
Let
$$
\Phi=
(\roman{Ad}(z_1)-\roman{Id},\dots,\roman{Ad}(z_n)-\roman{Id})
\colon
\g^n
@>>>
\g^n.
$$
Application of the functor $\roman {Hom}_{\Bobb R \pi}(\cdot,\g_{\phi})$
to (2.12) yields the cochain map
$$
\CD
\bold C(\Cal P, \g_{\phi}) \colon
@.
\g
@>{\delta_{\phi}^0}>>
\g^{2\ell}  \times \g^n
@>{\delta_{\phi}^1}>>
\g
\\
@.
@A{\roman{Id}}AA
@A{(\roman{Id},\Phi)}AA
@A{\roman{Id}}AA
\\
\bold C(\widetilde {\Cal P},\{\pi_j\}; \g_{\phi}) \colon
\quad
@.
\g
@>>{\delta_{\phi}^0}>
\g^{2\ell}  \times \g^n
@>>{\delta_{\phi}^1}>
\g
\endCD
$$
where the notation
$\delta_{\phi}^0$
and
$\delta_{\phi}^1$
is slightly abused.
It is obvious
that this chain map
factors through
$\bold C_{\roman{par}}(\Cal P, \g_{\phi})$.
A little thought reveals that this implies
the assertion. \qed
\enddemo

\proclaim{Proposition 4.5}
The tangent maps
$\roman T_e\alpha_{\phi}$
and $\roman T_{\phi}\rh$
make commutative the diagram
$$
\CD
@.
\roman T_eG
@>\roman T_e\alpha_{\phi}>>
\roman T_{\phi} \left (\roman{Hom}(F,G)_{\bold C}\right)
@>{\roman T_{\phi} \rh}>>
\roman T_{r (\phi)}G
\\
@.
@A{\roman{Id}}AA
@A{\roman R_{\phi}}AA
@A{\roman R_{r(\phi)}}AA
\\
\bold C_{\roman{par}}(\Cal P, \g_{\phi})\colon
@.
\g
@>>{\delta^0_{\phi}}>
\g^{2\ell} \times h_1 \times \dots\times h_n
@>>{\delta^1_{\phi}}>
\g,
\endCD
\tag4.5.1
$$
having its vertical arrows isomorphisms
of vector spaces.
\endproclaim

\demo{Proof}
In fact, the diagram (4.2)
restricts to the diagram (4.5.1). \qed
\enddemo

\smallskip
It is manifest that
the kernel of
the operator $\delta_{\phi}^1$
in $\bold C_{\roman{par}}(\Cal P, \g_{\phi})$
(cf. (4.4.1))
coincides with
the space
$Z_{\roman{par}}^1(\pi,\{\pi_j\};\g_{\phi})$
of parabolic 1-cocycles
with values in $\g_{\phi}$.

\beginsection 5. The extended moduli space

Let $\cdot$ be an invariant symmetric bilinear form on $\g$,
not necessarily positive definite or even nondegenerate.
The additional hypothesis of nondegeneracy
will be exploited only in Section 8 -- 10 below.
As in \cite\modus,
for a group $\Pi$, we denote by $(C_*(\Pi),\partial)$
the chain complex of its nonhomogeneous reduced normalized
bar resolution over the ground ring $R$.
When the relators
$z_1,\dots,z_n$ are added to
(2.1) we obtain the presentation
$$
\widehat {\Cal P} = \langle x_1,y_1,  \dots,x_\ell,y_\ell,z_1,\dots,z_n; r,
z_1,\dots,z_n
\rangle
\tag5.1
$$
of the fundamental group
$\widehat \pi =\pi_1(\widehat \Sigma)$
of the {\it closed\/}
surface $\widehat \Sigma$ resulting from capping of the $n$ boundaries.
Let $F$ be the free group on the generators
of $\Cal P$; notice that the generators of the latter
coincide with those of $\widehat {\Cal P}$.
We apply a variant of
the construction in \cite\modus\ to the presentation
$\widehat {\Cal P}$:
Let $O$ be the open $G$-invariant subset of the Lie algebra
$\g$ of $G$ where the exponential mapping is regular.
Define
the space
$\Cal H(\wide {\Cal P},G)_{\bold C}$
by means of the pull back square
$$
\CD
\Cal H(\wide {\Cal P},G)_{\bold C}
@>{(\widehat r,\overline z_1,\dots,\overline z_n)}>>
O \times C_1\times\dots\times C_n
\\
@V{\eta}VV
@VV{\roman{exp} \times \roman{Id} \times \dots\times \roman{Id}}V
\\
\roman{Hom}(F,G)_{\bold C}
@>>{(r,z_1,\dots,z_n)}> G\times C_1 \times \dots \times C_n,
\endCD
\tag5.2
$$
where
$\widehat r$ and $\overline z_1,\dots,\overline z_n$
denote the induced maps.
The space
$\Cal H(\wide {\Cal P},G)_{\bold C}$
is manifestly
a smooth manifold.
\smallskip
Let $c$ be an absolute 2-chain of $F$
which represents a 2-cycle for the group system
$(\pi;\pi_1,\dots,\pi_n)$.
Its image in the 2-chains
of the fundamental group
$\widehat \pi$
of the {\it closed\/} (!)
surface $\widehat \Sigma$ is then closed.
Write $\widehat\kappa \in \roman H_2(\widehat \pi)$ for its class.
When the genus $\ell$ is different from zero
$\widehat \pi$
is non-trivial and
the canonical map
from
$\roman H_2(\pi,\{\pi_j\})$
to $\roman H_2(\widehat \pi)$
is an isomorphism
identifying the fundamental classes.
Write
$$
E\colon F^2 \times \roman{Hom}(F,G) @>>> G^2
$$
for the evaluation map, and let
$$
\omega_c = \langle c, E^*\Omega\rangle,
$$
the result of pairing $c$ with the induced form,
cf.
\cite {\modus\  (13)}.
This is a $G$-invariant  2-form on
$\roman{Hom}(F,G)$.
In view of
\cite {\modus\  (15)}
we have
$$
d\omega_c = \langle \partial c, E^*\lambda\rangle.
\tag5.3
$$
We now apply a variant of the construction
in Theorem 1
of \cite \modus:
\smallskip
We pick
$c$ in such a way that
$$
\partial c = [r] - [z_1] - \dots - [z_n]
\tag5.4
$$
in the
chain complex $C_*(F)$ of the nonhomogeneous
reduced normalized
bar resolution
of $F$.
This can always be done,
cf. what is said about the chain (2.11)
at the end of Section 2 above.
In fact,
for $\ell \geq 1$,
the construction in \cite \modus,
applied to the presentation
$$
\langle x_1,y_1,  \dots,x_\ell,y_\ell; r^{\flat}\rangle,
\quad
r^{\flat} = [x_1,y_1] \dots [x_\ell,y_\ell] = rz_1^{-1} \dots z_n^{-1}
$$
of the fundamental group $\widehat \pi$
of $\widehat\Sigma$,
yields a 2-chain
$c_1$ with
$$
\partial c_1 = [rz_1^{-1} \dots z_n^{-1}]  \in C_1(F)
$$
Since
$$
\align
\partial [r\,z_n^{-1} \dots z_2^{-1}|z_1^{-1}]
&=
[z_1^{-1}]
- [r\,z_n^{-1} \dots z_1^{-1}]
+ [r\,z_n^{-1} \dots z_2^{-1}]
\\
\partial [r\,z_n^{-1} \dots z_3^{-1}|z_2^{-1}]
&=
[z_2^{-1}]
- [r\,z_n^{-1} \dots z_2^{-1}]
+ [r\,z_n^{-1} \dots z_3^{-1}]
\\
& \cdots
\\
\partial [r|z_n^{-1}]
&=
[z_n^{-1}]
- [r\,z_n^{-1}]
+ [r] ,
\endalign
$$
adding to $c_1$ the
2-chains coming into play on the
left-hand sides of these
equations, we arrive at a 2-chain $c_2$ with
$$
\partial c_2 = [r]+ [z_1^{-1}] + \dots + [z_n^{-1}]  \in C_1(F).
$$
Finally, subtracting the 2-chains $[z_1|z_1^{-1}],\dots,
[z_n|z_n^{-1}]$
we obtain the desired 2-chain $c$
satisfying (5.4) as asserted.
When $\ell=0$,
$c$ may be taken to be the negative of the sum of the chains
$$
[z_1\cdot \dots \cdot z_{n-1}|z_n],
\
[z_1\cdot \dots \cdot z_{n-2}|z_{n-1}],
\
\dots\ ,
\
[z_1|z_2].
$$
\smallskip
Henceforth we suppose that
the homotopy operator $h$ on the forms on $\g$
used to construct the various forms in
\cite\modus\ is the standard operator.
Then the
map $\psi$ from
$\g$ to $\g^*$, cf. Lemma 1 in \cite\modus,
boils down to the adjoint of the
symmetric bilinear form $\cdot$ on $\g$.
Let
$\beta = h(\roman{exp}^*(\lambda))$, and
define the  2-form
$\omega_{c,\wide {\Cal P}}$
on
$\Cal H(\wide {\Cal P},G)_{\bold C}$
by
$$
\omega_{c,\wide {\Cal P}}
=
\eta^*\omega_c
- \widehat r^* \beta.
\tag5.5
$$
Since $d\beta
=\roman{exp}^*(\lambda)$,
in view of (5.3) above,
$$
d\omega_{c,\wide {\Cal P}}
=
-\overline z_1^*\lambda
-
\dots
-
\overline z_n^*\lambda
\tag5.6
$$
where we do not distinguish in notation
between
$\lambda$ and its restrictions to the conjugacy classes
$C_1,\dots, C_n$.
Next,
let
$$
\mu = \psi \circ \widehat r \colon
\Cal H(\wide {\Cal P},G)_{\bold C}
@>>>
\g^*
$$
that is,
$\mu$ is the composite of
the map
$\widehat r$
from $\Cal H(\wide {\Cal P},G)_{\bold C}$ to $\g$
with the adjoint $\psi$ of the symmetric
bilinear form $\cdot$ from $\g$ to its dual;
here we do not distinguish in notation between a map into
$O$ and its composite with the inclusion into $\g$.
Recall
$$
\delta_G\lambda = - d \vartheta \quad \text{\cite{\modus\ (6)}} .
$$
As in \cite\modus, we write
$$
\mu^\sharp
\colon
\g
@>>>
C^{\infty}(\Cal H(\wide {\Cal P},G)_{\bold C})
$$
for the adjoint of $\mu$.
We now assert
$$
\delta_G
\omega_{c,\wide {\Cal P}}
=
d\mu^\sharp
-
\overline z_1^*\vartheta
-
\dots
-
\overline z_n^*\vartheta
\tag 5.7
$$
where we do not distinguish in notation
between
$\vartheta$ and its restrictions to the conjugacy classes
$C_1,\dots, C_n$.
Indeed,
cf. the proof of Theorem 2
of \cite \modus,
$$
\psi^{\sharp} = h \circ (\roman{exp}^* \vartheta - \delta_G(\beta))
$$
and
$$
\delta_G(\Omega) = \delta \vartheta \quad
\text{\cite{\modus\ (4)}}
$$
whence
$$
\align
\delta_G
\omega_{c,\wide {\Cal P}}
&=
\delta_G
(\eta^*\omega_c
- \widehat r^* \beta)
\\
&=
\eta^*\delta_G\omega_c
- \widehat r^* \delta_G\beta
\\
&=
\eta^*\langle \partial c, E^*\vartheta\rangle
- \widehat r^* \delta_G\beta
\\
&=
\widehat r^*(\roman{exp}^*\vartheta - \delta_G \beta)
-
\overline z_1^*\vartheta
-
\dots
-
\overline z_n^*\vartheta
\\
&=
d\mu^\sharp
-
\overline z_1^*\vartheta
-
\dots
-
\overline z_n^*\vartheta
\endalign
$$
as asserted.
The formulas (5.6) and (5.7)
show that
$\mu$ is somewhat like a momentum
mapping for the (non-closed)  2-form
$\omega_{c,\Cal P}$,
up to certain error terms.
In Section 7 below we shall
add appropriate forms
which will correct this error.
Before we can do so we need some preparation
to which the next Section is devoted.

\beginsection 6. A single conjugacy class

Let $C$ be a conjugacy class of $G$
and $\Cal O$  an adjoint orbit
which is mapped onto $C$ under the exponential
mapping from $\g$ to $G$.
Let $X,Y \in \g$.
The vector fields
$X_{\Cal O}$
and
$Y_{\Cal O}$
on
$\Cal O$ generated by $X$ and $Y$ are given by the assignment
to a point $Z\in \Cal O$ of
$[X,Z] \in \roman T_Z\Cal O$
and
$[Y,Z] \in \roman T_Z\Cal O$,
and the \lq\lq Kirillov\rq\rq\  form
$\omega$ on $\Cal O$ is given by the expression
$$
\omega_Z(X_{\Cal O}, Y_{\Cal O}) =
\omega_Z([X,Z],[Y,Z]) = [X,Y]\cdot Z = [Z,X]\cdot Y.
\tag6.1
$$
Notice that there is no need to assume the
symmetric bilinear form
$\cdot$ on $\g$ to be nondegenerate; just take
the  2-form $\omega$ on $\Cal O$ defined by (6.1).
For a point $p$ of $C$, an arbitrary tangent vector
is of the form
$$
Xp - pX =
(X- \roman {Ad}(p) X)p  \in \roman T_pC,
$$
where $ \cdot\, p$ and $p\, \cdot $
denote the effect of right and left translation
and where $X$ is an element of the Lie algebra $\g$, identified
with the tangent space $\roman T_eG$ of $G$ at $e$.
As before, let $\beta = h (\roman{exp}^*\lambda)$
where $\lambda$ denotes Cartan's fundamental  3-form on $G$.

\proclaim{Theorem 6.2}
The assignment
$$
\tau(Xp- pX, Yp - pY) =
\frac 12 (X \cdot \roman{Ad}(p)Y
-
Y \cdot \roman{Ad}(p)X),
\quad
p \in C,
\tag6.2.1
$$
yields an equivariant  2-form $\tau$ on $C$
having the property
$$
\roman{exp}^*\tau = \beta - \omega.
\tag6.2.2
$$
\endproclaim

\smallskip
Before proving the theorem, we spell out the following which will be crucial:

\proclaim{Corollary 6.3}
The  2-form $\tau$ satisfies the formulas
$$
\align
d\tau &= \lambda
\tag 6.3.1
\\
\delta_G\tau &= \vartheta \, .
\tag 6.3.2
\endalign
$$
\endproclaim

\demo{Proof of the Corollary}
Since $\omega$ is closed,
$$
\roman{exp}^*(d\tau)= d\roman{exp}^*\tau
=d(\beta-\omega) = d \beta = \roman{exp}^*\lambda
$$
whence
$d\tau = \lambda$.
Furthermore,
denote by $J$ the composite of the inclusion of $\Cal O$
into $\g$ with the adjoint of the given symmetric bilinear form;
formula
(6.1) says that
$$
\omega_Z(X_{\Cal O},Y_{\Cal O})
=
d(X \circ J)_Z(Y_{\Cal O}),
$$
that is,
with our definition
of the operator $\delta_G$
involving the negative (!) of the contraction operator,
cf. Section 1 of \cite\modus,
we have
$$
\delta_G\omega = - J^\sharp.
$$
On the other hand,
in view of the formula
$$
\delta_G\lambda = - d \vartheta
\quad
\text{\cite{\modus\ (4)}},
$$
on the whole Lie algebra $\g$, we get
$$
\align
\delta_G(\beta)
&=
\delta_G(h\roman{exp}^*\lambda)
\\
&=
-h(\delta_G(\roman{exp}^*\lambda))
\\
&=
-h(\roman{exp}^*\delta_G(\lambda))
\\
&=
h(\roman{exp}^* d(\vartheta))
\\
&=
hd(\roman{exp}^* (\vartheta))
\\
&=
\roman{exp}^* (\vartheta)
-dh(\roman{exp}^* (\vartheta))
\\
&=
\roman{exp}^* (\vartheta)
-\psi^\sharp.
\endalign
$$
Consequently, on $\Cal O$,
where $\psi$ amounts to $J$,
we obtain
$$
\roman{exp}^*(\delta_G\tau) =
\delta_G(\roman{exp}^*\tau) =
\delta_G(\beta-\omega)
=
\roman{exp}^* (\vartheta)
-\psi^\sharp
+ J^{\sharp}
=
\roman{exp}^* (\vartheta)
$$
whence
$$
\delta_G\tau = \vartheta
$$
holds on $C$
as asserted. \qed
\enddemo
\smallskip
\noindent{\smc Remark.}
Notice that the proof of (6.3) works whether or not
the restriction of the exponential mapping
to $\Cal O$
is a diffeomorphism.
\smallskip
We now begin with the preparations for the proof
of Theorem 6.2.
Recall \cite {\helgaboo~II.1.7} that the derivative
at $Z \in \g$
of the exponential mapping
$\roman{exp}$
from $\g$ to $G$ is given by the formula
$$
d \roman{exp}_Z
=
d(L_{\roman{exp} Z})_e
\circ
\frac{1-e^{-\roman{ad Z}}}{\roman{ad} Z}
\tag 6.4.1
$$
that is,
$$
d \roman{exp}_Z
=
d(L_{\roman{exp} Z})_e
\circ
\left(
1 - \frac 12
\roman{ad} Z
+ \frac 1{3!}
(\roman{ad} Z)^2
-
\frac 1{4!}
(\roman{ad} Z)^3
+\dots \right) .
\tag 6.4.2
$$
We now consider the exponential mapping
$\roman{exp}$ from $\Cal O$ to the corresponding conjugacy class $C$.
Let $p = \roman{exp} (Z) \in C \subseteq G$.
The above formula entails that the derivative
$$
d \roman{exp}_Z\colon
\roman T_Z\Cal O
@>>>
\roman T_pC
$$
of the exponential mapping
sends the tangent vector $[X,Z]$ to
the tangent vector
$$
Xp - pX =
(X-\roman {Ad}(p) X)p  \in \roman T_pC.
$$
Consequently
the statement of  Theorem 6.2 is equivalent to the following.

\proclaim{Lemma 6.5} The  2-form $\beta$ on $\Cal O$
is given by the formula
$$
\beta_Z([X,Z],[Y,Z])
=
[Z,X]\cdot Y
+
\frac 12
(X\cdot \roman{Ad}(p) Y -
Y\cdot \roman{Ad}(p) X).
$$
\endproclaim

Henceforth we write $[\cdot,\cdot,\cdot]$
for the triple product.

\demo{Proof}
As in \cite\modus,
write
$\rho = \roman{exp}^*\lambda$,
where $\lambda$ refers to the fundamental
3-form on $G$.
For simplicity, write $p_t =
\roman{exp}(tZ) \in G$.
We then have
$$
\align
2 \beta_Z([X,Z],[Y,Z])
&=
2 \int_0^1
\rho_{tZ}(Z,[X,tZ],[Y,tZ]) dt
\\
&=
\int_0^1
\left[p_t^{-1}(d\roman{exp}_{tZ}(Z)),
\roman{Ad}(p_t^{-1}) X - X,
\roman{Ad}(p_t^{-1}) Y - Y \right] dt
\\
&=
\int_0^1
\left[Z,
\roman{Ad}(p_t^{-1}) X - X,
\roman{Ad}(p_t^{-1}) Y - Y \right] dt
\\
&=
\int_0^1
\left[Z,X,Y\right] dt
+
\int_0^1
\left[Z,
\roman{Ad}(p_t^{-1}) X,
\roman{Ad}(p_t^{-1}) Y\right] dt
\\
&\quad
+
\int_0^1
[Z,\roman{Ad}(p^{-1}_t)X,-Y] dt
+
\int_0^1
[Z,-X,\roman{Ad}(p^{-1}_t)Y] dt
\\
&=
2[Z,X]\cdot  Y
\\
&\quad
+
\int_0^1
\left(
\roman{Ad}(\roman{exp}(tZ)X,Z] \cdot Y
-
[\roman{Ad}(\roman{exp}(tZ)Y,Z] \cdot X
\right)dt.
\endalign
$$
However,
$$
[\roman{Ad}(\roman{exp}(tZ)Y,Z]
=
-\left([Z,Y] + t [Z,[Z,Y]] + \frac {t^2}{2!} [Z,[Z,[Z,Y]]] + \dots
\right)
$$
whence
$$
\align
\int_0^1 [\roman{Ad}(\roman{exp}(tZ)Y,Z] dt
&=
-\left([Z,Y] + \frac 12 [Z,[Z,Y]] + \frac 1{3!} [Z,[Z,[Z,Y]]] + \dots
\right)
\\
&=
Y - e^{\roman {ad}(Z)}Y
\\
&=
Y - \roman {Ad}(p)Y
\endalign
$$
and likewise
$$
\int_0^1 [\roman{Ad}(\roman{exp}(tZ)X,Z] dt
=
X - \roman {Ad}(p)X
$$
whence the assertion. \qed
\enddemo

\beginsection 7. The completion of the construction

In view of (5.6), (5.7), and
(6.3) above, very little work remains to prove the following.

\proclaim{Theorem 7.1}
The equivariant  2-form
$$
\omega_{c,\wide {\Cal P},\bold C}
=
\omega_{c,\wide {\Cal P}}
+
\overline z_1^*\tau
+
\dots
+
\overline z_n^*\tau
\tag7.1.1
$$
on
$\Cal H(\wide {\Cal P},G)_{\bold C}$
is closed, and
the adjoint
$\mu^\sharp$ from
$\g$ to $C^{\infty}(\Cal H(\wide {\Cal P},G)_{\bold C})$ of the
smooth equivariant
map
$\mu$
from $\Cal H(\wide {\Cal P},G)_{\bold C}$ to $\g^*$
satisfies
the identity
$$
\delta_G\omega_{c,\wide {\Cal P},\bold C}
= d\mu^\sharp .
\tag7.1.2
$$
Consequently the difference
$\omega_{c,{\Cal P},{\bold C}}
-\mu^\sharp$
is an equivariantly closed
form
in
\linebreak
$(\Omega_G^{*,*}(\Cal H(\wide {\Cal P},G)_{\bold C});d,\delta_G)$
of total degree 2.
\endproclaim

The identity (7.1.2) says that,
for every $X \in \g$,
$$
-\omega_{c,{\Cal P},{\bold C}}(X_{\Cal H},\cdot\,)
= d (X \circ \mu),
$$
that is, $\mu$ is formally a momentum mapping
for the $G$-action on
$
\Cal H(\wide {\Cal P},G)_{\bold C}$,
with reference to
$\omega_{c,\wide {\Cal P},\bold C}$,
except that
the latter is not necessarily nondegenerate;
here we have written
$X_{\Cal H}$ for the vector field
on
$\Cal H(\wide {\Cal P},G)_{\bold C}$
induced by $X \in \g$.
\smallskip

\beginsection 8. Groupoids and the equivariantly
closed form

Let $\phi\in \roman{Hom}(F,G)_{\bold C}$
and suppose that $\phi(r)$ lies in the centre of $G$.
The above construction
(3.6), applied to the present data, yields an
alternating  bilinear form
$$
\omega_{\kappa,\cdot,\phi}\colon
\roman H^1_{\roman{par}}(\pi,\{\pi_j\};\g_{\phi})
\otimes
\roman H^1_{\roman{par}}(\pi,\{\pi_j\};\g_{\phi})
@>>>
\bold R
\tag8.1
$$
on
$\roman H^1_{\roman{par}}(\pi,\{\pi_j\};\g_{\phi})$
which is symplectic, i.~e. nondegenerate, provided that
$\cdot$ is nondegenerate.
Next, let
$$
\omega_{c,\bold C}
=
\omega_c
+
\overline z_1^*\tau
+
\dots
+
\overline z_n^*\tau .
\tag8.2
$$
This is a  2-form on
$\roman{Hom}(F,G)_{\bold C}$
whose restriction
to
$\roman{Hom}(\pi,G)_{\bold C}$
coincides with the restriction
of
$\omega_{c,\Cal P,\bold C}$
to
$\roman{Hom}(\pi,G)_{\bold C}$
where
$\roman{Hom}(\pi,G)_{\bold C}$
is viewed as a subspace of
$\Cal H(\wide {\Cal P},G)_{\bold C}$
as explained above.
Henceforth we denote by
$K_{\phi}$ the kernel of the derivative
$\roman T_{\phi}r$
occurring in {\rm (4.5.1)} above.
Our present goal is to prove the following.

\proclaim{Theorem 8.3}
Right translation identifies the restriction of the  2-form
$\omega_{c,\bold C}$ to
$K_{\phi}$
with the alternating bilinear form on
$Z^1_{\roman{par}}(\pi,\{\pi_j\};\g_{\phi})$ obtained as the composite
of
$\omega_{\kappa,\cdot,\phi}$
with the projection from
$Z^1_{\roman{par}}(\pi,\{\pi_j\};\g_{\phi})$
to
$\roman H^1_{\roman{par}}(\pi,\{\pi_j\};\g_{\phi})$.
\endproclaim

\proclaim{Key Lemma 8.4}
For an arbitrary real representation $V$ of $\pi$
with an invariant symmetric bilinear form
$\cdot$\,,
the
value of
the alternating bilinear form {\rm (3.6)}
on
$\roman H^1_{\roman{par}}(\pi,\{\pi_j\};V)$
for
two parabolic $V$-valued 1-cocycles
$u$ and $v$
with
$$
u(z_j) = z_j X_j - X_j
\quad
\text{and}
\quad v(z_j) = z_j Y_j - Y_j,
\qquad X_j, Y_j \in V,
\ 1 \leq j \leq n,
\tag8.4.1
$$
is given by the formula
$$
\omega_V([u],[v]) =
\langle c, u \cup v \rangle
+ \frac 12 \sum \left( X_j \cdot z_j Y_j  -Y_j \cdot z_j X_j \right).
\tag8.4.2
$$
\endproclaim

We postpone the proof for the moment
and now give the

\demo{Proof of Theorem 8.3}
The chosen 2-chain $c\in C_2(F)$ in
the nonhomogeneous reduced normalized bar resolution
for $F$
looks like
$$
c = \sum \nu_{j,k}[x_j|x_k].
$$
Define
the bilinear form
$\omega_{c,\cup,\phi}$
on
$$
Z^1(F,\g_{\phi}) =
C^1({\Cal P},\g_{\phi})
\quad (=C^1(\widehat {\Cal P},\g_{\phi}))
$$
(cf. (5.1))
by the explicit formula
$$
\omega_{c,\cup,\phi}(u,v) = \langle c, u \cup v \rangle
=\sum \nu_{j,k}
u(x_j) \cdot (\roman{Ad}(\phi(x_j))v(x_k)),
\quad
u,v \in Z^1(F,\g_{\phi}).
$$
By   \cite{\modus\ (4.6)},
the  2-form $\omega_c$ is the
right translation of the
antisymmetrization
of $\omega_{c,\cup,\phi}$
and hence
the two coincide
on
$K_{\phi}$.
This involves the
Alexander-Whitney diagonal map
\cite\maclaboo\ (VIII.9 Ex. 1, p. 248).
Inspection shows that,
still on $K_{\phi}$,
the remaining terms
$\overline z_1^*\tau,\cdots,\overline z_n^*\tau$
in (8.2)
correspond precisely to
the remaining
terms in (8.4.2),
with
$V = \g_{\phi}$. In view of Lemma 8.4,
this proves Theorem 8.3. \qed
\enddemo

The proof of the Key Lemma  relies on a detailed analysis of the
multiplicative structure of the cohomology of the group
system in terms of the fundamental groupoid
$\widetilde \pi= \Pi(\Sigma;p_0,p_1,\dots,p_n)$
introduced in Section 2 above.
We now explain this.
\smallskip
Let
$\widetilde F$ denote the  groupoid
which is free
on the generators of (2.4).
To have a neutral notation,
whenever necessary,
we shall write
$\widetilde \Pi$
for either
$\widetilde F$ or
$\widetilde \pi$;
accordingly we write
$\Pi$
for either
$F$ or
$\pi$.
As usual, view $G$ as a groupoid with a single object which we write $e$.
Write $\roman{Hom}(\widetilde \Pi,G)$
for the space of groupoid homomorphisms from
$\widetilde \Pi$ to $G$.
\smallskip
The assignments
$$
\align
i(e) &= p_0,
\quad
i(x_j) = x_j,
\quad
i(y_j) = y_j,
\quad
i(z_j) = \gamma_j a_j \gamma_j^{-1},
\\
\cor(p_j) &= e,
\quad
\cor(x_j) = x_j,
\quad
\cor(y_j) = y_j,
\quad
\cor(a_j) = z_j,
\quad
\cor(\gamma_j) =\roman{Id},
\endalign
$$
yield obvious functors
$i \colon \Pi \to \widetilde \Pi$
and
$\cor \colon \widetilde \Pi \to \Pi$
inducing a deformation retraction
of $\widetilde \Pi$ onto $\Pi$;
cf. e.~g. \cite\brownboo\ (6.5.13)
for this notion.
These functors induce
maps
$$
i^* \colon \roman{Hom}(\widetilde \Pi,G)
@>>>
\roman{Hom} (\Pi,G),
\quad
\cor^* \colon \roman{Hom}( \Pi,G)
@>>>
\roman{Hom} (\widetilde \Pi,G)
$$
which, for $\Pi = F$, are manifestly smooth.
We shall occasionally refer to
$i^*$ and
$\cor^*$ as
{\it restriction\/}
and
{\it corestriction\/},
respectively.
\smallskip
The obvious action of $G$
on $\roman{Hom}(\Pi,G)$
by conjugation
extends
to an action of
the group
$G^{\widetilde \Pi_0} \cong G\times \dots \times G$
($n+1$ copies of $G$)
on $\roman{Hom}(\widetilde \Pi,G)$
in the following way:
We denote by $s$ and $t$ the source and target mappings
from $\widetilde \Pi$ to $\widetilde \Pi_0$.
Given a homomorphism
$\alpha$ from
$\widetilde \Pi$ to $G$ and
$\vartheta \in G^{\widetilde \Pi_0}$,
the homomorphism
$\vartheta \alpha$ is defined by
$$
\vartheta \alpha(w) = \vartheta(t(w))\alpha(w)(\vartheta(s(w)))^{-1}.
$$
The orbit space for the
$G^{\widetilde \Pi_0}$-action on
$\roman{Hom}(\widetilde \Pi,G)$
will be denoted by
$\roman{Rep}(\widetilde \Pi,G)$.
\smallskip
As in the group case,
we denote by $\roman{Hom}(\widetilde \Pi,G)_{\bold C}$
the space of homomorphisms $\chi$ from $\widetilde \Pi$ to $G$
for which
the value $\chi(a_j)$ of each generator
$a_j$
lies in $C_j$, for $1 \leq j\leq n$.
The
$G^{\widetilde \Pi_0}$-action on
$\roman{Hom}(\widetilde \Pi,G)$
leaves
the subspace
$\roman{Hom}(\widetilde \Pi,G)_{\bold C}$
invariant, and we denote
by
$\roman{Rep}(\widetilde \Pi,G)_{\bold C}$
the orbit space for this action.
Since
$$
\chi(z_j) = \chi(\gamma_j a_j \gamma_j^{-1})
=
\chi(\gamma_j) \chi(a_j) \chi(\gamma_j)^{-1},\quad
j = 1,\dots, n,
$$
the condition \lq $\chi (a_j) \in C_j$\rq\
is equivalent to
the condition
\lq $\chi (z_j) \in C_j$\rq,
which we used in the group case.

\proclaim{Proposition 8.5}
The restriction
mapping
induces bijections
$$
i^* \colon \roman{Rep}(\widetilde \Pi,G)
@>>>
\roman{Rep} (\Pi,G),
\quad
i^* \colon \roman{Rep}(\widetilde \Pi,G)_{\bold C}
@>>>
\roman{Rep} (\Pi,G)_{\bold C}.
$$
\endproclaim

Thus we can study the structure of
$\roman{Rep} (\Pi,G)_{\bold C}$
by looking at
$\roman{Rep}(\widetilde \Pi,G)_{\bold C}$
instead.
In particular,
this remark applies
to the infinitesimal structure,
in the following way:
Write $\beta$ for the nonhomogeneous unreduced normalized
bar resolution.
The
retraction $\cor$ from
$\widetilde \Pi$ to $\Pi$
induces
a deformation retraction
from
the nerve
$N\widetilde \Pi$
of
$\widetilde \Pi$
to
the nerve
$
N\Pi$
of $\Pi$,
and a canonical section for the latter
is of course
induced by the injection $i$
from
$\Pi$ to $\widetilde \Pi$;
hence
$\cor$ induces
a deformation retraction
from
$\beta
\widetilde \Pi$
to
$\beta
\Pi$.
In particular,
$\beta \widetilde \Pi$
yields a free resolution
of $R$
in the category of left
$R\Pi$-modules
of the kind
written $A$
in Section 1 above.
In fact,
write $\widetilde \pi_{\partial}$ for the free subgroupoid
of $\widetilde F$
having  $p_1,\dots,p_n$
as objects
and
$a_1,\dots,a_n$
as free generators
for its morphisms;
this groupoid may also be viewed as a subgroupoid
of $\widetilde \pi$,
and we do not distinguish in notation
between the two subgroupoids.
Abstractly,
$\widetilde \pi_{\partial}$
amounts
of course to a disjoint union
of the $n$ free cyclic groups
$\pi_1,\dots,\pi_n$,
and
$$
\beta \widetilde \pi_{\partial}
=
\oplus_{j=1}^n \beta \pi_j
\subseteq
\beta \widetilde \pi
$$
in such a way that extension of scalars yields an injection
of
$B=\oplus_{j=1}^n R\pi \otimes_{\pi_j}\beta \pi_j$
of $R\pi$-complexes
onto
a direct summand of
$\beta\widetilde \pi$.
The
$R\pi$-complex $B$
plays exactly the same role
as that denoted in
Section 1 above
by the same symbol, we have the split exact sequence
(1.1) at our disposal, and
the quotient
$$
\beta(\widetilde \pi,
\widetilde \pi_{\partial})
=
A \big / B
=
\beta\widetilde \pi
\big / B
$$
computes the relative cohomology
$\roman H^*(\widetilde \pi,
\widetilde \pi_{\partial}; \cdot)$.
With the present interpretation
of $(A,B)$ as resolution
over
the group system
$(\pi;\pi_1,\dots,\pi_n)$,
the relative cohomology
$\roman H^*(\widetilde \pi,
\widetilde \pi_{\partial}; \cdot)$
{\it coincides\/} with the cohomology
$\roman H^*(\pi,
\{\pi_j\}; \cdot)$
of the group system
$(\pi;\pi_1,\dots,\pi_n)$,
though.
In particular,
the standard formula for
the
diagonal in the nonhomogeneous unreduced normalized
bar resolution
$\beta \widetilde \pi$
for $\widetilde \pi$
yields
the multiplicative structure
of
$\roman H^*(\pi,\{\pi_j\}; \cdot)$.

\demo{Proof of the Key Lemma 8.4}
Let
$u$ and $v$
be parabolic $V$-valued 1-cocycles
on $\beta \pi$ so that
(8.4.1) holds;
the calculation of the value
$\omega_V([u],[v])$ of
the pairing (3.6) and hence (8.1)
may now be split into the following steps:
\roster
\item
{\bf Extension:}
The composites
$u'= u \circ \cor$ and $v'= v \circ \cor$
of $u$ and $v$, respectively,
with the retraction $\cor$
from
$\beta \widetilde \pi$
to
$\beta \pi$
yields
extensions to
parabolic $V$-valued
1-cocycles $u'$ and $v'$
on $\beta \widetilde \pi$.
\item
{\bf Normalization:}
Normalize $u'$ and $v'$
to obtain
groupoid cocycles $\widetilde u$ and $\widetilde v$
which take the value zero on the peripheral part $B$
of the resolution, cf. Section 1 above;
notice that this amounts to the requirement that
$\widetilde u$ and $\widetilde v$
vanish on the $n$ boundary circles $S_1,\dots, S_n$.
\item
{\bf Lifting:}
Lift the 2-chain $c \in C_2(F)$
of the nonhomogeneous reduced normalized bar resolution
for $F$
to a 2-chain $\widetilde c \in C_2(\widetilde F)$ of
the nonhomogeneous reduced normalized bar resolution
for $\widetilde F$
which
(i) passes to a relative cycle
for $(\widetilde \pi, \widetilde \pi_{\partial})$
and which (ii)
under
the retraction $\cor$
from
$C_*(\widetilde F)$ to $C_*(F)$
goes to $c$.
\item
{\bf Computation:}
The value
$\omega_V([u],[v])$
is then computed by the formula
$$
\omega_V([u],[v])
=
\langle \widetilde c, \widetilde u \cup \widetilde v \rangle.
\tag8.6
$$
\endroster
We note that the groupoid description
is crucial for the normalization in step 2;
such a normalization would be impossible
for the group cocycles $u$ and $v$, and
a calculation
of the value
$\omega_V([u],[v])$
directly in terms of $u$ and $v$
would lead to a mess.
Moreover,
by a purely formal
reasoning in
the relative cohomology of the pair
$(A,B)$ or, what amounts to the same,
of the pair $(\Sigma, \partial \Sigma)$, the right-hand
side of (8.6) is well defined
and yields the pairing (3.6).
\smallskip
There is no more need to comment on step 1,
and we now explain the other steps.
\smallskip
{\bf Step 2.}
With reference to (8.4.1),
let $X$ and $Y$ be the groupoid 0-cocycles
defined by
$$
X(p_0) = 0,
\quad
X(p_j) = X_j,
\quad
Y(p_0) = 0,
\quad
Y(p_j) = Y_j,
\quad
1 \leq j \leq n.
$$
Define $\widetilde u$ and $\widetilde v$ by
$$
\widetilde u = u' - \delta X = u \circ \cor - \delta X,
\quad
\widetilde v = v' - \delta Y = v \circ \cor - \delta Y.
$$
\smallskip
{\bf Step 3.}
View $c$ as a 2-chain of
$\widetilde F$
by the embedding of
$F$ into $\widetilde F$
and let
$$
\widetilde c =
c + \sum_j\left([\gamma_j^{-1}|\gamma_j a_j] -
[\gamma_ja_j|\gamma_j^{-1}]
\right) .
\tag8.7
$$
Then
$$
\partial \widetilde c =
[\widetilde r] - [a_1] - \dots - [a_n]
$$
whence, in particular, $\widetilde c$ is
manifestly
a relative 2-cycle
for $(\widetilde \pi, \widetilde \pi_{\partial})$.
Notice that
$c$ itself is {\it not\/}
a relative cycle
for $(\widetilde \pi, \widetilde \pi_{\partial})$.
\smallskip
{\bf Step 4.} By definition
$$
\align
\omega_V([u],[v])
&=
\langle \widetilde c, \widetilde u \cup \widetilde v \rangle =
\\
&
\phantom{+}
\langle c, (u'-\delta X) \cup (v'-\delta Y)\rangle
\tag8.8
\\
&
+
\sum_j \langle [\gamma_j^{-1}|\gamma_j a_j],
(u'-\delta X) \cup (v'-\delta Y)\rangle
\tag8.9
\\
&
-
\sum_j
\langle [\gamma_ja_j|\gamma_j^{-1}],(u'-\delta X) \cup (v'-\delta Y)\rangle
\tag8.10
\endalign
$$
The term (8.8)
is the sum
$$
\langle c, u' \cup v'\rangle
-
\langle c, u' \cup \delta Y\rangle
-
\langle c, \delta X \cup v'\rangle
+
\langle c, \delta X \cup \delta Y\rangle
\tag8.11
$$
Since $c$ is a group chain,
$
\langle c, u' \cup v'\rangle
$
equals
$\langle c, u \cup v\rangle$.
Further, since $u'$ is a cocycle,
$$
\align
\langle c, u' \cup \delta Y\rangle
&=
\langle c, -\delta(u' \cup Y)\rangle
=
\langle -\partial c, u' \cup Y\rangle
\\
&=
\sum_j   \langle [z_j], u' \cup Y\rangle
=
\sum_j   u(z_j) \cdot Y (p_0) = 0
\endalign
$$
since
$Y (p_0) = 0$.
Likewise,
the last two terms in
(8.11) involve $X(p_0)$, which is zero as well,
so we are left with
$\langle c, u \cup v\rangle$
as the value of the entire sum (8.11)
which computes (8.8).
The term (8.9) involves $n$ factors which may be computed as
$$
\langle [\gamma_j^{-1}|\gamma_j a_j],
(u'-\delta X) \cup (v'-\delta Y)\rangle
=
(u'-\delta X)(\gamma_j^{-1}) \cdot
(\gamma_j^{-1} ((v'-\delta Y) (\gamma_j a_j))).
$$
However,
$u'(\gamma_j^{-1}) = u(\cor \gamma_j^{-1})
= u(e)$ which is zero,
$$
\align
\delta X(\gamma_j^{-1}) &= X(\partial \gamma_j^{-1}) = X(p_0 -p_j) = -X_j,
\\
v'(\gamma_j a_j) &= v(\cor(\gamma_j a_j)) = v(z_j) = z_j Y_j - Y_j,
\\
\delta Y (\gamma_j a_j) &= Y(\partial (\gamma_j a_j)) = Y(a_j p_j - p_0)
=z_j Y_j,
\endalign
$$
and
$\gamma_j$
and $\gamma_j^{-1}$ act as the identity
on $V$.
Hence (8.9) equals
$$
\sum_j X_j \cdot(z_j Y_j - Y_j - z_j Y_j) =
-\sum_j X_j\cdot Y_j.
$$
Finally the term (8.10) involves $n$ factors which may be computed as
$$
\align
\langle [\gamma_ja_j|\gamma_j^{-1}],&(u'-\delta X) \cup (v'-\delta Y)\rangle
\\
&=
\left(
(u'-\delta X)
(\gamma_ja_j)
\right)
\cdot
(\gamma_j a_j)\left((v'-\delta Y)(\gamma_j^{-1})\right)
\\
&=
\left(
u(\cor(\gamma_j a_j)) - \delta X(\gamma_j a_j)
\right)
\cdot
(\gamma_j a_j)
\left(
v(\cor (\gamma_j^{-1})) - \delta Y(\gamma_j^{-1})
\right)
\\
&=
\left(
u(z_j) - X(a_j p_j - p_0)
\right)
\cdot
(\gamma_j a_j)
\left(
- Y(\partial \gamma_j^{-1})
\right)
\\
&=
\left(
z_j X_j - X_j - z_jX_j
\right)
\cdot
z_j
Y_j
\\
&= - X_j \cdot z_j Y_j
\endalign
$$
Consequently
$$
\omega_V([u],[v])
=
\langle c, u \cup v \rangle
+ \sum X_j \cdot (z_j Y_j - Y_j).
$$
Likewise,
$$
\omega_V([v],[u])
=
\langle c, v \cup u \rangle
+ \sum Y_j \cdot (z_j X_j - X_j).
$$
By antisymmetry,
$$
\align
2\omega_V([u],[v])
&=
\omega_V([u],[v])
-
\omega_V([v],[u])
\\
&
=
2\langle c, u \cup v \rangle
+ \sum X_j \cdot (z_j Y_j - Y_j)
-
\sum Y_j \cdot (z_j X_j - X_j)
\\
&=
2\langle c, u \cup v \rangle
+ \sum \left( X_j \cdot z_j Y_j  -Y_j \cdot z_j X_j \right)
\endalign
$$
This proves the key Lemma.\qed
\enddemo

\beginsection 9. Stratified symplectic structures

Suppose that is $G$ compact and that the symmetric
bilinear form $\cdot$\ on $\g$
is nondegenerate.
Recall that the notion of a stratified symplectic space
has been introduced in \cite\sjamlerm.

\proclaim {Theorem 9.1}
With respect to the decomposition according to $G$-orbit types,
the space
$\roman{Rep}(\pi,G)_{\bold C}$
inherits a structure of a stratified symplectic space.
\endproclaim

In fact, the argument
for
the main result of
\cite\sjamlerm\
shows that
each connected
component of a reduced space of the kind considered
inherits a structure of a stratified symplectic space.
In the setting
of \cite\sjamlerm\
the hypothesis of properness
is used {\it only \/}
to
guarantee that the reduced space is in fact connected.
In our situation, we know a priori
that the reduced space is connected.

\proclaim{Corollary 9.2}
The space
$\roman{Rep}(\pi,G)_{\bold C}$
has a unique open, connected, and dense stratum.
\endproclaim

In fact, this follows at once from \cite\sjamlerm\ (5.9).
The stratum mentioned in the corollary is called the {\it top stratum\/}.
Thus there is a certain subgroup $T$ of $G$, unique up to conjugacy,
such that every $\phi \in
\roman{Hom}(\pi,G)_{\bold C}$
representing a point of the top stratum
has stabilizer $Z_{\phi}$ conjugate to $T$.
In many cases, $T$ is just the centre of $G$,
and the top stratum consists of representations
which are {\it irreducible\/}
in the sense that the stabilizer has
Lie algebra the Lie algebra of the centre of $G$
(but the top stratum
may be smaller than the space of irreducible representations).
See \cite\singulat\ for details.
There may be {\it no\/}
irreducible representations at all,
though.
This happens for example when $\pi$ is abelian
and $G$ non-abelian.

\beginsection 10. Relationship with gauge theory constructions

For $G$ compact, and closed \cite{\atibottw,\gurupone}\ as well as punctured
\cite{\bisguron,\bisgurtw}\ surfaces, the symplectic structure
on the top stratum
of the corresponding space of
gauge equivalence classes of
flat connections
has been described using methods of gauge theory.
In this section we
show that the usual identification
of this space
with the corresponding space
of representations
identifies
the symplectic structures on the top strata.
This extends what is done
in Section 6 of \cite\smooth\  where, for the case of a closed surface,
the symplectic structures on all strata
have been shown to correspond to each other.
\smallskip
Let $G$ be a general, not necessarily compact Lie group.
The compactness hypothesis
will be made at the appropriate stage.
Our  surface $\Sigma$ is compact, with $n\geq 0$ boundary circles
$S_1,\dots, S_n$
and chosen base point $p_0$.
Write $\Sigma^\bullet$ for the corresponding
{\it punctured\/} surface
with base point $p_0$
which contains
$\Sigma$
as a based deformation retract
in such a way that each $S_j$ is a circle about the corresponding puncture.
We
do not exclude the case
of a closed surface
and we agree that in this case
$\Sigma^\bullet$
and $\Sigma$ coincide.
\smallskip
Let $\xi \colon P \to
\Sigma^\bullet$
be a
flat
principal $G$-bundle,
having the structure group $G$ act from the right as usual,
and pick a base point $\widehat p_0$ of $P$ with
$\xi(\widehat p_0) = p_0$.
In many cases,
for example when $\pi$ is a free group or when $G$ is simply connected,
$\xi$ will be topologically trivial.
Write $\Cal A(\xi)$ for the space of connections on $\xi$.
The assignment to a gauge transformation
$\nu$ on $\xi$ of $x_\nu \in G$
defined by $\nu(\widehat p_0) = \widehat p_0 x_\nu$
furnishes a homomorphism from the group
$\Cal G(\xi)$ of gauge transformations onto $G$.
Among the various descriptions of
the space
$\Omega^*(\Sigma^\bullet,\roman{ad}(\xi))$
of forms with values in
the adjoint bundle
$$
\roman{ad}(\xi)
\colon P \times _G \g
@>>> \Sigma^\bullet
$$
we shall take that in terms of $G$-invariant horizontal
$\g$-valued forms on $P$;
we note that here $G$ acts on
its Lie algebra $\g$
from the {\it left\/}
by the adjoint representation
$\roman{Ad}\colon G \to \roman{Aut}(\g)$
as usual.
\smallskip
For
a smooth closed path
$w \colon I \to \Sigma^\bullet$
defined on the unit interval $I$, with starting point
$p_0 \in \Sigma^\bullet$,
the {\it holonomy\/}
$\roman{Hol}_{w,\widehat p_0}(A) \in G$
of $A$
{\it with reference to\/}
$\widehat p_0$
is defined by
$$
\widehat w (1)
= \widehat p_0\, \roman{Hol}_{w,\widehat p_0}(A) \in P
$$
where $\widehat w$ refers to the {\it horizontal\/}
lift of $w$ having starting  point $\widehat p_0$.
For $b \in G$, we denote by
$L_b$ the operation of left translation
from $\g$ to $\roman T_bG$.
\smallskip
We maintain the notation of Sections 2 and 4 above.
In particular,
$F$ is the free group on the generators
$x_1,y_1,\dots,x_{\ell}, y_{\ell},z_1,\dots, z_n$
of (2.1) and, with an abuse
of notation,
the corresponding closed (edge) paths in $\Sigma$
representing these generators
are denoted by the same symbols.
The assignment to
a connection $A$ of
the point
$$
\left(\roman{Hol}_{x_1,\widehat p_0}(A),
\roman{Hol}_{y_1,\widehat p_0}(A),
\dots,
\roman{Hol}_{x_\ell,\widehat p_0}(A),
\roman{Hol}_{y_\ell,\widehat p_0}(A),
\roman{Hol}_{z_1,\widehat p_0}(A), \dots,
\roman{Hol}_{z_n,\widehat p_0}(A)\right)
$$
of $G^{2\ell + n}$
yields a smooth  map
$$
\rho
\colon
\Cal A(\xi)
@>>>
G^{2\ell + n} =\roman{Hom}(F,G)
$$
which is $\Cal G(\xi)$-equivariant in the sense that
$$
\rho(\nu A)= x_\nu \rho(A) x_\nu^{-1},
\quad
\text{
for every gauge transformation}\quad \nu;
$$
this map is
referred to as
{\it Wilson loop mapping\/}
in \cite\smooth, where a comment is made as to the
appropriate interpretation of the
property of $\rho$ being smooth.
The  restriction of $\rho$ to the subspace
$\Cal F (\xi)$ of flat connections
yields the standard map from
$\Cal F (\xi)$
to $\roman{Hom}(\pi,G)$,
viewed as a subspace of
$\roman{Hom}(F,G)$ via the projection from $F$ to $\pi$,
and this map
depends only on the choice of $\widehat p_0$
but
{\it not\/} on the choices of closed paths
representing the generators of $F$.
The induced map from
the space
of gauge equivalence
classes of flat connections on $\xi$ to
(the corresponding open and closed subset of)
the representation space
$\roman{Rep}(\pi,G)=\roman{Hom}(\pi,G)\big/ G$
is then independent of the choice of
$\widehat p_0$.
\smallskip
Now let $A$ be a flat connection, and let
$
\phi = \rho(A)
$
be the corresponding homomorphism from $\pi$ to $G$.
Write
$\Omega^1=\Omega^1(\Sigma^\bullet,\roman{ad}(\xi))$
and $\roman T_{\phi}=\roman T_{\phi} G^{2\ell + n}$.
Consider  the universal cover
$\widetilde \Sigma^\bullet$ of
$\Sigma^\bullet$, and suppose things arranged in such a way that
$\pi$ acts on
$\widetilde \Sigma^\bullet$
and $\widetilde\Sigma$
from
the {\it right\/}. As usual, this action is related
to the corresponding action of $\pi$ from the left
by
$x\widetilde p = \widetilde p x^{-1}$,
for
$x \in \pi$ and $\widetilde p \in \widetilde\Sigma^\bullet$.
After a choice $o$ of base point
of $\widetilde \Sigma^\bullet$
over $p_0$ has been made,
there is a canonical smooth $\pi$-equivariant map
$\sigma$ from
$\widetilde \Sigma^\bullet$
to $P$ over the identity mapping of
$\Sigma^\bullet$;
here $\pi$ acts on  $P$ via $\phi$
from the right.
Explicitly, given a point
$\widetilde p$  of $\widetilde \Sigma^\bullet$,
let $\widetilde w$ be a smooth path in
$\widetilde\Sigma^\bullet$
from
$o$
to
$\widetilde p$,
write $w$ for its projection into
$\Sigma^\bullet$,
and let $\widehat w$ be the horizontal lift
of $w$, with reference to $A$ and
$\widehat p_0$;
then
$\sigma(\widetilde p)$ equals the end point of  $\widehat w$.
In particular,
$\sigma$
induces an isomorphism
of principal bundles
from
$
\xi_\phi \colon
\widetilde \Sigma^\bullet \times _\phi G
\to
\Sigma^\bullet
$
to $\xi$
identifying
the obvious flat connection $A_\phi$
on
$\xi_\phi$
with $A$, but this fact will not be needed below;
here
$\widetilde \Sigma^\bullet \times _\phi G$
is the space arising
from $\widetilde \Sigma^\bullet \times G$
by identifying points of the kind
$(p,\phi(x) b)$
and
$(p\phi(x), b)$ as usual.

\smallskip
In Section 4 above,
the Lie algebra $\g$  has been viewed as a {\it left\/}
$\pi$-module  via $\phi$,
written $\g_{\phi}$ and,
via right translation,
we identified the tangent space
$\roman T_{\phi}$
with
the space
$C^1(\Cal P, \g_{\phi})$
of $\g_{\phi}$-valued {\it left\/} 1-cochains on (2.2).
We shall use the same notation $\g_{\phi}$ for $\g$, viewed
as
a {\it right\/}
$\pi$-module.
There is {\it no\/} conflict of notation
since the left and right $\pi$-actions on $\g$ are related by
$$
xX= X x^{-1} = \roman{Ad}(\phi(x))X,
\quad X \in \g,\ x \in \pi.
$$
Accordingly, there are two notions of
$\g_{\phi}$-valued
1-cocycles:
A {\it left\/} 1-cocycle
is a function $u$ from $\pi$ to $\g$
satisfying
$u(xy) = u(x) + x u(y)\ (=u(x) + \roman{Ad}(\phi(x))u(y))$
while a {\it right\/} 1-cocycle
is a function $v$ from $\pi$ to $\g$
satisfying
$v(xy) = (v(x))y + v(y)$ where $(v(x))y=\roman{Ad}(\phi(y)^{-1})(v(x))$.
Likewise,
via $\phi$, the group $\pi$
acts  on  the de Rham complex
$(\Omega^*(\widetilde\Sigma^\bullet, \g_\phi),d)$
from the right,
and we can take
invariants $(\Omega^*(\widetilde\Sigma^\bullet, \g_\phi),d)^{\pi}$.
The operator
of covariant derivative
$d_A$
is a differential on
$\Omega^*(\Sigma^\bullet, \roman{ad}(\xi))$,
and the above map $\sigma$ induces an isomorphism
$$
(\Omega^*(\Sigma^\bullet, \roman{ad}(\xi)),d_A)
@>>>
(\Omega^*(\widetilde\Sigma^\bullet, \g_\phi),d)^{\pi}
\tag10.1
$$
where on the right-hand side
the symbol $d$ refers to the usual de Rham coboundary operator.

\proclaim{Proposition 10.2}
For a flat connection $A$, with $\phi= \rho(A)$, under the identification
of
$\roman T_A(\Cal A(\xi)) =
\Omega^1$
with
$\Omega^1(\widetilde\Sigma^\bullet, \g_\phi)^{\pi}$
via {\rm (10.1)} (in degree 1)
and of
$\roman T_\phi$ with $\g^{2\ell +n}= C^1(\Cal P, \g_{\phi})$
via right translation,
the derivative of $\rho$ at $A$ assigns to a closed
$\pi$-invariant
$\g_\phi$-valued 1-form $\vartheta$
on
$\widetilde\Sigma^\bullet$
the $\g_{\phi}$-valued
left 1-cocycle $u_{\vartheta}$ for $\pi$
given by the formula
$$
u_\vartheta (x)
=
\int_{xo}^{o} \vartheta,
\quad
\text{for}\  x \in \pi;
$$
here the integral is taken along any smooth path
in $\widetilde\Sigma^\bullet$
from
$xo$ to $o$,
and
$\widetilde\Sigma^\bullet$
is viewed
as a left $\pi$-space.
This assignment induces the usual isomorphism from
$\roman H^1_A(\Sigma^\bullet, \roman{ad}(\xi))$
onto
$\roman H^1(\pi, \g_\phi)$.
\endproclaim

\demo{Proof}
Theorem 2.7 of \cite\smooth\  entails that,
at an arbitrary
connection $A$,
not necessarily flat,
with
$$
\rho(A) =
(a_1,b_1,\dots, a_{\ell}, b_{\ell}, c_1,\dots,c_n) \in G^{2\ell + n},
$$
the differential
$
d\rho(A)
\colon
\roman T_A\Cal A(\xi)
@>>>
\roman T_{\rho(A)} G^{2\ell + n}
$
of $\rho$ is given by the assignment to
$\vartheta \in \Omega^1(\Sigma^\bullet,\roman{ad}(\xi))= \roman T_A\Cal A(\xi)$
of
the vector
$$
\left(
L_{a_1} \int_{\widehat x_1} \vartheta,
L_{b_1} \int_{\widehat y_1} \vartheta,
\dots,
L_{a_\ell} \int_{\widehat x_\ell} \vartheta,
L_{b_\ell} \int_{\widehat y_\ell} \vartheta,
L_{c_1} \int_{\widehat z_1} \vartheta,
\dots,
L_{c_n} \int_{\widehat z_n} \vartheta
\right)
$$
in
$
\roman T_{a_1} G
\times
\roman T_{b_1} G
\times
\dots
\times
\roman T_{a_\ell} G
\times
\roman T_{b_\ell} G
\times
\roman T_{c_1} G
\times
\dots
\times
\roman T_{c_n} G
$.
However,
when $A$ is flat,
the tangent map of $\rho$,
combined with
the inverse of (10.1)
(in degree 1) and
with left translation
from
$\roman T_{\phi}$
to
$\g^{2\ell + n}$,
looks like
$$
\Omega^1(\widetilde\Sigma^\bullet, \g_\phi)^{\pi}
@>>>
\Omega^1
@>>>
\roman T_{\phi}
@>>>
\g^{2\ell + n};
$$
it is given by the assignment to
a $\pi$-invariant
$\g_\phi$-valued 1-form $\vartheta$
on
$\widetilde\Sigma^\bullet$
of
$$
v_\vartheta=
\left(
\int_{\widehat x_1} \vartheta,
\int_{\widehat y_1} \vartheta,
\dots,
\int_{\widehat x_\ell} \vartheta,
\int_{\widehat y_\ell} \vartheta,
\int_{\widehat z_1} \vartheta,
\dots,
\int_{\widehat z_n} \vartheta
\right)
\in
\g^{2\ell + n}
$$
where,
with an abuse of notation,
$\widehat x_j$,
$\widehat y_j$,
and
$\widehat z_k$
refer to the unique lifts
in
$\widetilde\Sigma^\bullet$
with reference to $o$
of,
respectively,
the closed paths
$x_j$,
$y_j$,
and
$z_k$ in $\Sigma$.
This assignment
is in fact the degree one twisted integration mapping
from
$(\Omega^*(\widetilde\Sigma^\bullet, \g_\phi),d)^{\pi}$
to
the cellular cochains
$(C^*(\Sigma, \g_\phi),d)$
with {\it local coefficients\/}
determined by $\phi$,
cf. Section 4 of \cite\smooth\
and what is said below.
In particular,
the cellular 1-cocycles
$Z^1(\Sigma, \g_\phi)$
with local coefficients
coincide with the
$\g_\phi$-valued
right 1-cocycles
for $\pi$.
Thus for a closed 1-form $\vartheta$,
the cochain $v_\vartheta$
yields a
$\g_\phi$-valued
right 1-cocycle for $\pi$.
Since integration of a closed 1-form
on a simply connected space
does not depend on the choice of path
but only on the endpoints,
this 1-cocycle
assigns to $x \in \pi$ the integral
$$
v_\vartheta (x) = \int_{o}^{o x} \vartheta
$$
taken along any smooth path from
$o$ to $o x$.
\smallskip
The tangent map of $\rho$,
combined with
the inverse of (10.1)
(in degree 1)
and
{\it right\/} translation
from
$\roman T_{\phi}$
to
$\g^{2\ell + n}$,
is given by the assignment to
a
$\pi$-invariant
$\g_\phi$-valued 1-form $\vartheta$
on
$\widetilde\Sigma^\bullet$
of the vector
$$
\left(
\roman{Ad}(\phi(x_1))\int_{\widehat x_1} \vartheta,
\dots,
\roman{Ad}(\phi(y_\ell))\int_{\widehat y_\ell} \vartheta,
\roman{Ad}(\phi(z_1))\int_{\widehat z_1} \vartheta,
\dots,
\roman{Ad}(\phi(z_n))\int_{\widehat z_n} \vartheta
\right)
$$
in
$\g^{2\ell + n}$.
For a closed 1-form $\vartheta$,
this yields the
$\g_{\phi}$-valued
left 1-cocycle
$u_\vartheta$
for $\pi$
assigning
the value
$$
u_\vartheta (x)
=
\roman{Ad}(\phi(x)) v_\vartheta (x)
=
\roman{Ad}(\phi(x))
\int_{o}^{o x} \vartheta
=
\roman{Ad}(\phi(x))
\int_{o}^{x^{-1}o} \vartheta
=
\int_{xo}^{o} \vartheta
$$
to $x \in \pi$. \qed
\enddemo

Let $A$ be a flat connection on $\xi$, and let $\phi= \rho(A)$.
The
right-hand side of (10.1)
involves only the vector space of de Rham forms and
the homomorphism $\phi$ but no longer  the connection $A$
explicitly.
Thus we can completely do away with
$\xi$ and the flat connection $A$ and work entirely in terms
of $\phi$ and the local system it defines in the following way,
where for the sake
of clarity we proceed in somewhat greater generality
than actually needed:
Let $V$ be a real representation of $\pi$;
in the application below, $V$ will be $\g_\phi$.
The de Rham complex
$(\Omega^*(\widetilde\Sigma^\bullet, V),d)$ inherits an obvious action of $\pi$
(whether or not we take the left or right incarnation thereof
will not matter any more since both lead to the same result),
and we can consider its $\pi$-invariants
$(\Omega^*(\widetilde\Sigma^\bullet, V),d)^{\pi}$;
by means of an isomorphism of the kind (10.1),
the
cohomology
$\roman H^*_{\roman{equiv}}(\widetilde\Sigma^\bullet,V)$
of
$(\Omega^*(\widetilde\Sigma^\bullet, V),d)^{\pi}$
actually computes the cohomology
of $\Sigma^\bullet$  with values in the corresponding flat vector bundle
as usual but this is not important here. Integration
$$
(\Omega^*(\widetilde\Sigma^\bullet, V),d)^{\pi}
@>>>
(C_{\roman{local}}^*(\Sigma^\bullet, V),d)
$$
into the cochains
$$
(C_{\roman{local}}^*(\Sigma^\bullet, V),d)
=
(C^*(\widetilde\Sigma^\bullet, V),d)^{\pi}
$$
on $\Sigma^\bullet$ with local coefficients determined by $V$
assigns, in particular, to a $\pi$-invariant
$V$-valued closed 1-form $\vartheta$
on
$\widetilde\Sigma^\bullet$
the $V$-valued
left 1-cocycle $u_{\vartheta}$ for $\pi$
given by the formula
$$
u_\vartheta (x)
=
\int_{xo}^{o} \vartheta,
\quad
\text{for}\  x \in \pi,
$$
and this association induces the standard isomorphism
from
$\roman H^1_{\roman{local}}(\Sigma^\bullet,V)$
onto
$\roman H^1(\pi,V)$.

\proclaim{Proposition 10.3}
When
$\vartheta$
is  a compactly supported
closed $\pi$-invariant
$V$-valued 1-form
on
$\widetilde\Sigma^\bullet$,
$u_{\vartheta}$
is a parabolic 1-cocycle, that is,
for $1 \leq j \leq n$, there is $X_j \in V$ such that
$$
u_{\vartheta}(z_j) = z_j X_j - X_j.
$$
\endproclaim

\demo{Proof}
Let $1 \leq j\leq n$,
and pick a point $s_j$ of
$\widetilde \Sigma^{\bullet}$.
Then
$$
\align
u_\vartheta (z_j)
&=
\int_{z_jo}^{o} \vartheta
=
\int_{z_jo}^{z_j s_j} \vartheta
+
\int_{z_js_j}^{s_j} \vartheta
+
\int_{s_j}^{o} \vartheta
\\
&=
z_j X_j - X_j +
\int_{z_js_j}^{s_j} \vartheta
\endalign
$$
where
$X_j =\int_o^{s_j} \vartheta$.
However,
since $\vartheta$ is compactly supported
it vanishes in a neighborhood of the
punctures  whence, for a suitable choice
of $s_j$,
the integral
$\int_{z_js_j}^{s_j} \vartheta$ is zero. \qed \enddemo

\smallskip\noindent
{\smc Remark 10.4.}
By pushing the boundary circles
further towards the punctures if necessary,
in (10.3) above,
 we can  in fact assume that
$\vartheta$ vanishes on
$\Sigma^\bullet \setminus \Sigma$
and hence in particular on
the boundary circles $S_1,\dots, S_n$.
For $1 \leq j \leq n$,
the point $s_j$ may then be taken to be
the end point of the unique lift
$\widehat \gamma_j$
of $\gamma_j$
having starting point $o$
so that $s_j$ is a pre-image of $p_j$
and the integral
$X_j =\int_o^{s_j} \vartheta$
may be taken along
$\widehat \gamma_j$
where the notation in Section 2 above is in force.
More generally,
given
{\it finitely many\/}
compactly supported forms
we may still assume that
things have been arranged in such a way that these forms
vanish
on
$\Sigma^\bullet \setminus \Sigma$.
\smallskip

Suppose $V$ endowed with a $\pi$-invariant
symmetric bilinear form $\cdot$\ ;
together with the wedge product
of forms it induces a bilinear pairing
$$
\wedge
\colon\Omega^1(\widetilde\Sigma^\bullet, V)^{\pi}
\otimes
\Omega^1(\widetilde\Sigma^\bullet, V)^{\pi}
@>>>
\Omega^2(\Sigma^\bullet,\Bobb R).
$$
We can now spell out the main technical statement of the present section.

\proclaim{Theorem 10.5}
For compactly supported
closed $\pi$-invariant
$V$-valued 1-forms
$\eta$ and $\vartheta$
on
$\widetilde\Sigma^\bullet$,
$$
\int_{-\Sigma^\bullet} \eta \wedge \vartheta
=
\omega_V([u_\eta],[u_\vartheta])
\tag10.5.1
$$
where
$\omega_V$
is the skew-symmetric bilinear pairing {\rm (3.6)}.
\endproclaim

Here
$-\Sigma^\bullet$
refers to
$\Sigma^\bullet$
endowed with the orientation opposite to that determined by
the disk $D$ in Section 2 above.

\demo{Proof}
Write
$\widetilde u_\eta$ and
$\widetilde u_\vartheta$
for the corresponding groupoid cocycles
arising from
$u_\eta$ and
$u_\vartheta$
by the normalization procedure
in step 2 of the proof of Lemma 8.4.
In view of (8.6),
$$
\omega_V([u_\eta],[u_\vartheta])
=
\langle \widetilde c, \widetilde u_\eta \cup \widetilde u_\vartheta \rangle
$$
where $\widetilde c$ is
the groupoid chain arising from $c$
by the construction
in step 3 of the proof of Lemma 8.4.
\smallskip
We may suppose that $\eta$ and $\vartheta$
vanish
on
$\Sigma^\bullet \setminus \Sigma$, cf. Remark 10.4 above.
Then
$$
\int_{\Sigma^\bullet} \eta \wedge \vartheta
=
\int_{\Sigma} \eta \wedge \vartheta .
$$
The cell decomposition of
$\Sigma$ induces a cell decomposition
of its universal cover
$\widetilde \Sigma$.
Extending earlier notation,
we denote by $\widehat x_j$,
$\widehat y_j$,
$\widehat a_k$,
and
$\widehat \gamma_k$
the unique lifts
in
$\widetilde\Sigma$
with reference to $o$
of,
respectively,
the edge paths
$x_j$,
$y_j$,
$a_k$,
and
$\gamma_k$
in $\Sigma$.
These
edge paths, together with their
left translates
under the $\pi$-action,
constitute the 1-cells of the cell decomposition
of
$\widetilde\Sigma$.
Inspection shows that
$$
\gathered
\widetilde u_\eta(x_j) = - \int_{\widehat x_j} \eta,
\quad
\widetilde u_\eta(y_j) = - \int_{\widehat y_j} \eta,
\quad
\widetilde u_\vartheta(x_j) = - \int_{\widehat x_j} \vartheta,
\quad
\widetilde u_\vartheta(y_j) = - \int_{\widehat y_j} \vartheta,
\\
\widetilde u_\eta(\gamma_k) = - \int_{\widehat \gamma_k} \eta,
\quad
\widetilde u_\eta(a_k) =0,
\quad
\widetilde u_\vartheta(\gamma_k) = - \int_{\widehat \gamma_k} \vartheta,
\quad
\widetilde u_\vartheta(a_k) =0.
\endgathered
$$
In other words,
the groupoid cocycles
$\widetilde u_\eta$
and
$\widetilde u_\vartheta$
coincide precisely with
the
$\pi$-equivariant $V$-valued cellular 1-cocycles
$\widehat u_\eta$
and
$\widehat u_\vartheta$
arising from $\eta$ and $\vartheta$, respectively,
under the integration mapping
$$
(\Omega^*(\widetilde\Sigma, V),d)
@>>>
(C_{\roman{cell}}^*(\widetilde\Sigma, V),d)
$$
from $V$-valued
de Rham forms to $V$-valued cellular cochains
$\widetilde\Sigma$,
perhaps up to a sign
depending
on how things have been adjusted
but irrelevant for us since the formula
(10.5.1) does not depend on this sign;
it depends on the orientation of $\Sigma$, though.
\smallskip
The disk $D$ mentioned in Section 2 above lifts to a disk
$\widehat D$
in $\widetilde\Sigma$, and the left translates
of $\widehat D$ under $\pi$ constitute the 2-cells of
$\widetilde\Sigma$.
Furthermore,
$$
\int_{\Sigma} \eta \wedge \vartheta
=
\int_{\widehat D} \eta \wedge \vartheta
$$
where on the right-hand side
the
wedge product
$\eta \wedge \vartheta$
is
viewed as a 2-form on
$\widetilde \Sigma^\bullet$.
\smallskip
The cellular chains of $\widetilde \Sigma$
are given by (2.5).
Comparing (2.7)
with (8.7),
viewing $\widetilde c$
as a groupoid cochain for $\widetilde \pi$,
 and
exploiting the standard fact that
the integration mapping
from de Rham cohomology
to usual
(cellular or singular) cohomology
is compatible with multiplicative structures,
we conclude that
$$
\int_{\widehat D} \eta \wedge \vartheta
=
-\langle \widetilde c, \widetilde u_\eta \cup \widetilde u_\vartheta \rangle.
$$
This completes the proof.
A closer look shows that,
after a suitable
triangulation of
$\Sigma$
arising from the
nerve of $\widetilde \pi$,
a suitably chosen 2-chain $\widetilde c$
amounts to the negative
of the corresponding subdivision of
the cellular chain $D$.
For related matters see for example what is said in \cite\maclaboo\
(IV.5 p. 119) and on p. 495 of \cite\eilmac. \qed
\enddemo
\smallskip
Now suppose $G$ compact and connected.
Then the Wilson loop mapping $\rho$
identifies the moduli space
$N(\xi)$ of
gauge equivalence classes of flat connections
on $\xi$ with
(an open and closed subset of)
the representation space
$\roman{Rep}(\pi,G)=\roman{Hom}(\pi,G)\big/ G$.
Write
$N(\xi)_{\bold C}$ for the pre-image of
the subspace
$\roman{Rep}(\pi,G)_{\bold C}=\roman{Hom}(\pi,G)_{\bold C}\big/ G$
of
$\roman{Rep}(\pi,G)=\roman{Hom}(\pi,G)\big/ G$
under this identification;
thus
$N(\xi)_{\bold C}$
consists of gauge equivalence classes of
flat connections $A$
so that, for $1 \leq j \leq n$, the
holonomy along some small circle
about
the $j$'th puncture
lies in the chosen conjugacy class
$C_j$.
\smallskip
Suppose $\cdot$ nondegenerate. For a point
$[\phi]$
of the top stratum
$\roman{Rep}(\pi,G)_{\bold C}^{\roman{top}}$
of $\roman{Rep}(\pi,G)_{\bold C}$
(cf. (9.2) above),
in view of (4.4) and (4.5),
a choice of representative
$\phi$
induces an isomorphism
$\lambda_\phi$
from
$\roman H^1_{\roman{par}}(\pi,\{\pi_j\};\g_\phi)$
onto the (usual smooth) tangent space
$\roman T_{[\phi]}(\roman{Rep}(\pi,G)_{\bold C}^{\roman{top}})$;
this isomorphism
is independent of the choice
of $\phi$ in the sense that,
for every $x \in G$,
the composite
$$
\roman H^1_{\roman{par}}(\pi,\{\pi_j\};\g_\phi)
@>{\roman{Ad}(x)}>>
\roman H^1_{\roman{par}}(\pi,\{\pi_j\};\g_{x\phi})
@>{\lambda_{x\phi}}>>
\roman T_{[\phi]}(\roman{Rep}(\pi,G)_{\bold C}^{\roman{top}})
$$
coincides with $\lambda_\phi$.
This makes precise the folklore statement that
\lq the tangent space is the first cohomology group with coefficients
in the corresponding Lie algebra representation\rq;
details for the special case with no
punctures have been worked out in
Section 7 of \cite\smooth.
\smallskip
Under the Wilson loop mapping,
the top stratum
$\roman{Rep}(\pi,G)_{\bold C}^{\roman{top}}$
of $\roman{Rep}(\pi,G)_{\bold C}$
corresponds to the
subspace
$N(\xi)^{\roman{top}}_{\bold C}$
of
points
$[A]$ of
$N(\xi)_{\bold C}$
whose representatives $A$ have minimal stabilizer
subgroup (in the group of gauge transformations).
For a flat connection $A$, write
$\roman H^1_{A,c}(\Sigma^\bullet,\roman{ad}(\xi))$
for the subgroup
of
the
first
cohomology group
$\roman H^1_A(\Sigma^\bullet,\roman{ad}(\xi))$
generated by classes of
compactly supported
1-forms.
A choice of representative $A$
of a point $[A]$ of
$N(\xi)^{\roman{top}}_{\bold C}$
induces an isomorphism
$\lambda_A$
from
$\roman H^1_{A,c}(\Sigma^\bullet,\roman{ad}(\xi))$
onto
the
(usual smooth)
tangent space
$\roman T_{[A]}(N(\xi)^{\roman{top}}_{\bold C})$
and, for a gauge transformation $\nu$,
the composite
$$
\roman H^1_{A,c}(\Sigma^\bullet,\roman{ad}(\xi))
@>{\nu}>>
\roman H^1_{\nu A,c}(\Sigma^\bullet,\roman{ad}(\xi))
@>{\lambda_{\nu A}}>>
\roman T_{[A]}(N(\xi)^{\roman{top}}_{\bold C})
$$
coincides with $\lambda_A$.
\smallskip
Let $A$
be a flat connection on $\xi$ representing a
point of
$N(\xi)_{\bold C}$
and let $\phi = \rho(A)$.
The isomorphism
(10.1) identifies
$\roman H^1_{A,c}(\Sigma^\bullet,\roman{ad}(\xi))$
with
the subgroup
$\roman H^1_{\roman{equiv},c}(\widetilde\Sigma^\bullet,\g_\phi)$
of
$\roman H^1_{\roman{equiv}}(\widetilde\Sigma^\bullet,\g_\phi)$
generated by classes of compactly supported closed equivariant 1-forms.
In view of (10.3),
integration identifies
$\roman H^1_{\roman{equiv},c}(\widetilde\Sigma^\bullet,\g_\phi)$
with
$\roman H^1_{\roman{par}}(\pi,\{\pi_j\};\g_\phi)$.
Suppose in addition that $[A]$
lies in
$N(\xi)^{\roman{top}}_{\bold C}$.
Up to signe, the gauge theory description of
the 2-form on
$N(\xi)^{\roman{top}}_{\bold C}$
induced by
$\cdot$
is given
on
the tangent space
$\roman T_{[A]}(N(\xi)^{\roman{top}}_{\bold C})$
by the left-hand side of
(10.5.1) with
$V = \g_\phi$.
On the tangent space
$\roman T_{\phi}\roman{Rep}(\pi,G)_{\bold C}^{\roman{top}}$,
via right translation,
the 2-form
induced by
$\cdot$
is given
on $\roman H^1_{\roman{par}}(\pi,\{\pi_j\};\g_\phi)$
by the right-hand side of
(10.5.1) with
$V = \g_\phi$.
Theorem 10.5
implies at once that
the two forms correspond, up to sign.
This identifies the gauge theory description
of the symplectic form
with the
representation space
description
given in the present paper.
Notice that the given identification is independent
of the symplecticity.
\smallskip
A similar statement can be made
at an arbitrary point $[A]$ of
$N(\xi)_{\bold C}$
and the corresponding point $[\phi]$ of
$\roman{Rep}(\pi,G)_{\bold C}$
where $\phi = \rho(A)$.
Choices of $\phi$ and $A$ still determine
linear maps
of the kind $\lambda_A$ and $\lambda_\phi$
but these maps
will in general no longer be isomorphisms.
More precisely, $\lambda_\phi$
induces an isomorphism of the subspace
$\roman H^1_{\roman{par}}(\pi,\{\pi_j\};\g_\phi)^{Z_\phi}$
of invariants
onto the smooth tangent space at $[\phi]$
of the stratum in which
$[\phi]$ lies
where $Z_\phi \subseteq G$
refers to the stabilizer of $\phi$;
a corresponding statement
can be made for $A$.
This has been worked out for the closed case
(no punctures) in Section 7 of \cite\smooth.
\smallskip
When $G$ is not compact,
while the statement of (10.5) is still available,
there is no good space of gauge equivalence
classes of flat connections
nor is there
a good space of representations
since there are orbits which are not closed.
The appropriate generalization of the present results
should involve certain categorical or algebraic quotients.

\bigskip
\centerline{References}
\medskip\noindent
\widestnumber\key{999}

\ref \no \atiyboo
\by M. Atiyah
\book The geometry and physics of knots
\publ Cambridge University Press
\publaddr Cambridge, U. K.
\yr 1990
\endref

\ref \no  \atibottw
\by M. Atiyah and R. Bott
\paper The Yang-Mills equations over Riemann surfaces
\jour Phil. Trans. R. Soc. London  A
\vol 308
\yr 1982
\pages  523--615
\endref

\ref \no \biecktwo
\by R. Bieri and B. Eckmann
\paper Relative homology and Poincar\'e  duality
for group pairs
\jour J. Pure and Applied Algebra
\vol 13
\yr 1978
\pages 277--319
\endref

\ref \no \bisguron
\by I. Biswas and K. Guruprasad
\paper Principal bundles on open surfaces and invariant functions
on Lie groups
\jour Int. J. of Math.
\vol 4
\yr 1993
\pages 535--544
\endref

\ref \no \bisgurtw
\by I. Biswas and K. Guruprasad
\paper On some geometric invariants associated to the space of flat
connections
\jour Can. J. of Math.  (to appear)
\endref

\ref \no \bottone
\by R. Bott
\paper On the Chern-Weil homomorphism and the continuous cohomology of
Lie groups
\jour Advances
\vol 11
\yr 1973
\pages  289--303
\endref

\ref \no \botshust
\by R. Bott, H. Shulman, and J. Stasheff
\paper On the de Rham theory of certain classifying spaces
\jour Advances
\vol 20
\yr 1976
\pages 43--56
\endref

\ref \no \brownboo
\by R. Brown
\book Elements of modern topology
\publ McGraw-Hill
\publaddr London
\yr 1968
\endref

\ref \no \eilmac
\by S. Eilenberg and S. Mac Lane
\paper Relations between homology and homotopy groups of spaces
\jour Ann. of Math.
\vol 46
\yr 1945
\pages 480--509
\endref

\ref \no \fockrotw
\by V. V. Fock and A. A. Rosly
\paper Flat connections and polyubles
\paperinfo in: Modern problems in quantum field theory, strings,
and quantum gravity, Kiev 1992
\jour Teor. Mat. Fiz.
\vol 95
\yr 1993
\pages 228--238
\endref

\ref \no \goldmone
\by W. M. Goldman
\paper The symplectic nature of fundamental groups of surfaces
\jour Advances
\vol 54
\yr 1984
\pages 200--225
\endref

\ref \no \gurupone
\by K. Guruprasad
\paper Flat connections, geometric invariants, and the symplectic nature
of the fundamental group of surfaces
\jour Pacific J. of Math.
\vol 162
\yr 1984
\pages 45--55
\endref

\ref \no \gururaja
\by K. Guruprasad and C. S. Rajan
\paper Group cohomology and the symplectic structure
on the moduli space of representations
\paperinfo preprint, McGill University, 1995
\endref

\ref \no \helgaboo
\by S. Helgason
\book Differential geometry, Lie groups, and symmetric spaces
\bookinfo Pure and Applied Mathematics, vol. 80
\publ Academic Press
\publaddr New York $\cdot$ San Francisco $\cdot$ London
\yr 1978
\endref

\ref \no \modus
\by J. Huebschmann
\paper Symplectic and Poisson structures of certain moduli spaces
\jour Duke Math. J. (to appear)
\paperinfo hep-th 9312112
\endref

\ref \no \modustwo
\by J. Huebschmann
\paper Symplectic and Poisson structures of certain moduli spaces. II.
Projective
representations of cocompact planar discrete groups
\jour Duke Math. J. (to appear)
\paperinfo dg-ga/9412003
\endref

\ref \no \singula
\by J. Huebschmann
\paper The singularities of Yang-Mills connections
for bundles on a surface. I. The local model
\jour Math. Z., to appear \finalinfo dg-ga/9411006
\endref
\ref \no \singulat
\by J. Huebschmann
\paper The singularities of Yang-Mills connections
for bundles on a surface. II. The stratification
\jour Math. Z., to appear \finalinfo dg-ga/9411007
\endref
\ref \no \smooth
\by J. Huebschmann
\paper
Smooth structures on
moduli spaces of central Yang-Mills connections
for bundles on a surface
\paperinfo preprint 1992,
 dg-ga/9411008
\endref

\ref \no \poisson
\by J. Huebschmann
\paper
Poisson
structures on certain
moduli spaces
for bundles on a surface
\jour Annales de l'Institut Fourier
\vol 45
\yr 1995
\pages 65--91
\endref
\ref \no \locpois
\by J. Huebschmann
\paper Poisson geometry of flat connections
for {\rm SU(2)}-bundles on surfaces
\jour Math. Z., to appear \finalinfo hep-th/9312113
\endref

\ref \no \srni
\by J. Huebschmann
\paper
Poisson geometry of certain
moduli spaces
\paperinfo
Lectures delivered at the 14th winter school,
Czech Republic, Srni, January 1994
\jour Rendiconti del circolo matematico di Palermo,
to appear
\endref

\ref \no \huebjeff
\by J. Huebschmann and L. Jeffrey
\paper Group cohomology construction of symplectic forms
on certain moduli spaces
\jour Int. Math. Research Notices
\vol 6
\yr 1994
\pages 245--249
\endref

\ref \no \jeffrone
\by L. Jeffrey
\paper
Extended moduli spaces of flat connections
on Riemann surfaces
\jour Math. Ann.
\vol 298
\yr 1994
\pages 667--692
\endref

\ref \no \jeffrtwo
\by L. Jeffrey
\paper Symplectic forms on moduli spaces
of flat connections on 2-manifolds
\paperinfo preprint alg-geom/9411008, to appear in {\it Proceedings
of the Georgia International Topology Conference}, Athens, Ga.
1993, ed. by W. Kazez
\endref

\ref \no \jeffrthr
\by L. Jeffrey
\paper
Group cohomology construction of the cohomology of
moduli spaces
of flat connections
on 2-manifolds
\jour Duke Math. J.
\vol 77
\yr 1995
\pages 407--429
\endref

\ref \no \karshone
\by Y. Karshon
\paper
An algebraic proof for the symplectic
structure of moduli space
\jour Proc. Amer. Math. Soc.
\vol 116
\yr 1992
\pages 591--605
\endref

\ref \no \maclaboo
\by S. Mac Lane
\book Homology
\bookinfo Die Grundlehren der mathematischen Wissenschaften
 No. 114
\publ Springer
\publaddr Berlin $\cdot$ G\"ottingen $\cdot$ Heidelberg
\yr 1963
\endref

\ref \no \narasesh
\by M. S. Narasimhan and C. S. Seshadri
\paper Stable and unitary vector bundles on a compact Riemann surface
\jour Ann. of Math.
\vol 82
\yr 1965
\pages  540--567
\endref

\ref \no \shulmone
\by H. B. Shulman
\book Characteristic classes and foliations
\bookinfo Ph. D. Thesis
\publ University of California
\yr 1972
\endref

\ref \no \sjamlerm
\by R. Sjamaar and E. Lerman
\paper Stratified symplectic spaces and reduction
\jour Ann. of Math.
\vol 134
\yr 1991
\pages 375--422
\endref
\ref \no \trottone
\by H. F. Trotter
\paper Homology of group systems with applications to knot theory
\jour Annals of Math.
\vol 76
\yr 1962
\pages 464--498
\endref
\ref \no  \weiltwo
\by A. Weil
\paper  Remarks on the cohomology of groups
\jour
Ann. of Math.
\vol 80
\yr 1964
\pages  149--157
\endref

\ref \no \weinstwe
\by A. Weinstein
\paper The symplectic structure on moduli space
\paperinfo  The A. Floer memorial volume,
H. Hofer, C. Taubes, A. Weinstein, and E. Zehnder, eds.,
Birkh\"auser Verlag,
Basel $\cdot$ Boston $\cdot$ Berlin, 1995,
pp. 627--635
\endref

\enddocument